\documentclass{article}
\usepackage[dvips]{graphics}

\usepackage{amssymb}
\usepackage{amsmath}
\usepackage{amsthm}
\usepackage{amsbsy}
\usepackage[dvips]{graphics}
\usepackage{amscd}

\theoremstyle{plain}

\newtheorem{theo}{Theorem}[section]
\newtheorem{pr}[theo]{Proposition} 
\newtheorem{lem}[theo]{Lemma}
 
\newtheorem{de}[theo]{Definition}
\newtheorem{de-pr}[theo]{Definition-Proposition}

\theoremstyle{remark}
\theoremstyle{example}
\newtheorem{remark}[theo]{Remark}
\newtheorem{example}[theo]{Example}
\def\R{\mathbb{R}}

\def\N{\mathbb{N}}
\def\C{\mathbb{C}}
\def\Z{\mathbb{Z}}
\def\Q{\mathbb{Q}}

\begin{document}

\title{An extension of the Burau representation to a mapping class group
  associated to Thompson's group $T$}

\author{
\begin{tabular}{cc}
 Christophe Kapoudjian &  Vlad Sergiescu\\
\small \em Laboratoire Emile Picard, UMR 5580\ 
&\small \em Institut Fourier BP 74, UMR 5582\\
\small \em University of Toulouse
III &\small \em University of Grenoble I\\
\small \em 118, route de Narbonne &\small \em \\
\small \em 31062 Toulouse cedex 4, France 
&\small \em 38402 Saint-Martin-d'H\`eres cedex, France\\
\small \em e-mail: {\tt ckapoudj@picard.ups-tlse.fr}
& \small \em e-mail: {\tt Vlad.Sergiescu@ujf-grenoble.fr} \\
\end{tabular}
}

\date{October 2002}
\maketitle

\begin{center}
{\bf Abstract}
\end{center}
We study some aspects of the geometric
representation theory of the Thompson and Neretin groups, suggested by their analogies with the diffeomorphism groups of the circle. We prove that the Burau representation of the Artin braid groups extends to a
mapping class group $A_T$ related to Thompson's group $T$ by a short exact
sequence $B_{\infty}\hookrightarrow A_T\rightarrow T$, where $B_{\infty}$ is
the infinite braid group. This {\it non-commutative} extension abelianises to
a central extension $0\rightarrow \Z\rightarrow
A_T/[B_{\infty},B_{\infty}]\rightarrow T\rightarrow 1$ detecting the
 {\it discrete} version $\overline{gv}$ of the Bott-Virasoro-Godbillon-Vey class. A morphism
 from the above non-commutative extension to a reduced Pressley-Segal
 extension is then constructed, and the class $\overline{gv}$ is realised as a pull-back of the reduced Pressley-Segal
 class. A similar program is carried out for an extension of the Neretin group
 related to the {\it combinatorial} version of the Bott-Virasoro-Godbillon-Vey
 class.

\section{Introduction}

\hspace{0.4cm} The purpose of this work is to study some aspects of the geometric
representation theory of the Thompson and Neretin groups. Strikingly enough,
this turns out to be linked with the Burau representation of the classical braid
groups.\\

Recall that Thompson's group $T$ (cf. \cite{ca}) is the group of piecewise linear dyadic
homeomorphisms of the circle; it is finitely presented and simple. Another way
to look at its elements is as germs near the boundary of {\it planar partial
automorphisms} of the binary tree (a {\it planar partial automorphism} is
defined outside a finite subtree and preserves the cyclic order at each
vertex). As for Neretin's group $N$ (cf. \cite{ne}), it is defined as the the group of germs
of {\it all} the partial automorphisms of the binary tree; it is an
uncountable simple group (cf. \cite{ka2}).\\
Both groups $T$ and $N$ have been observed to present analogies with the
group $\mbox{Diff}(S^1)$ of orientation-preserving diffeomorphisms of the circle. Thus, $T$ looks,
up to some extent, as a ``lattice" in $\mbox{Diff}(S^1)$, while $N$ as a
$p$-adic (and combinatorial) analogue.\\

A remarkable point of this analogy concerns the Bott-Virasoro-Godbillon-Vey
class. This is a two dimensional (continuous) cohomology class $gv$ belonging to
$H^2(\mbox{Diff}(S^1);\R)$, whose derivative -- a cohomology class of the Lie
algebra $\mbox{Vect}(S^1)$ of vector fields on the circle -- was introduced by
Gelfand and Fuks; it corresponds to the Virasoro algebra, the universal
central extension of $\mbox{Vect}(S^1)$. A relative class $\overline{gv}$ for the group $T$ was
introduced and studied in \cite{gh-se} by E. Ghys
and the second author. Note that it is a
$\Z$-valued cohomology class. Moreover, a $\Z/2\Z$-valued class
$\overline{GV}$ has been defined for $N$ by the first author (\cite{ka3}).\\

A basic aspect of the Godbillon-Vey class is its relation with representation
theory. Thus, Pressley and Segal (\cite{pr-se}, Chapter 6) introduced representations of $\mbox{Diff}(S^1)$ in the restricted linear group $GL_{res}(L^2(S^1))$ (where
$L^2(S^1)$ is the Hilbert space of square integrable functions of the
circle). The class $e^{i.gv}$ is then essentially the pull-back of a class in
$H^2(GL_{res}(L^2(S^1));\C^*)$, which can be described as follows: there is an
extension
$$(NC)_{PS} \mbox{\hspace{1cm}} 1\rightarrow {\mathfrak T}\longrightarrow
\tilde{\mathfrak E}
\longrightarrow GL_{res}(L^2(S^1)) \rightarrow 1,$$ where ${\mathfrak T}$ is
the group of determinant operators on the Hardy subspace $L^2(S^1)_+$. After
dividing by the kernel of the determinant morphism from ${\mathfrak T}$ to
$\C^*$, one obtains a central extension detecting the desired cohomology
class.\\
By the machinery of the second quantization, the group $GL_{res}({\cal H})$ of
a Hilbert space ${\cal H}$ admits a projective representation in the fermionic
Fock space $\Lambda({\cal H})$: the spinor representation. Composed with the
Pressley-Segal representations, it provides projective representations of
$\mbox{Diff}(S^1)$. At the infinitesimal level, one obtains unitarisable
heighest weight modules of the Virasoro algebra (see \cite{ne"}
Chapter 7).\\

The first aim of this work is to construct similar representations for $T$
and $N$, inducing the classes $\overline{gv}$ and $\overline{GV}$
respectively. These classes have non-commutative versions introduced in
\cite{gr-se} and \cite{ka3} by P. Greenberg and the present authors; these are extensions
$$ (NC)_T \mbox{\hspace{1cm}} 1\rightarrow B_{\infty} \longrightarrow A_T  \longrightarrow T\rightarrow
1, $$
$$ (NC)_N \mbox{\hspace{1cm}} 1\rightarrow {\mathfrak S}_{\infty} \longrightarrow A_N  \longrightarrow N\rightarrow
1.$$
When divided by the commutators of $B_{\infty}$ and ${\mathfrak
  S}_{\infty}$, they produce central extensions which correspond
to $\frac{\overline{gv}}{2}$ and $\overline{GV}$ respectively.\\
A further aim of this paper is to
investigate the link of the above {\it non-commutative} extensions with the
Pressley-Segal extension $(NC)_{PS}$. By {\it non-commutative extensions} we mean
extensions of groups which provide central extensions after an abelianisation
process. Before stating our main result, we emphasize that such non-commutative
extensions appeared in various contexts (see \cite{gr-se2}, \cite{mo}, \cite{ts} and
\cite{ka4}), and are related to interesting central extensions:\\

1. In \cite{gr-se2}, the authors study an extension of $\widetilde{PSL_2(\R)}$ (the universal cover of $PSL_2(\R)$) by a group of piecewise
  $\widetilde{PSL_2(\R)}$ homeomorphisms of $\R$. By
  abelianisation, they get the Steinberg extension from $K$-theory.\\
\indent 2. Tsuboi (\cite{ts}) considers an extension of $\mbox{Diff}(S^1)$ by
  a group of area preserving homeomorphisms of the unit disk $D^2$. He proves
  that after abelianisation one obtains the Euler class of $\mbox{Diff}(S^1)$.\\
\indent 3. In \cite{ka4} one describes an action of the Neretin group $N$ on a tower of
  moduli spaces of real stable curves. Lifting this action to the universal
  cover provides a non-commutative extension of $N$ by the
  fundamental group of the tower. Its abelianisation is a non-trivial central extension of $N$ by
  $\Z/2\Z$.\\
\indent 4. Morita (\cite{mo}) looks at the modular group $\Gamma_g$ as an
  extension of the Torelli goup by the Siegel modular group $Sp_{2g}(\Z)$. The
  (non-central!) extension obtained by abelianisation captures the
  tautological classes of Morita and Mumford.\\ 

We stress that all these extensions have a geometrically defined middle term
as well.\\

The construction of $(NC)_T$, as performed in \cite{gr-se}, amounts to an
embedding $T\rightarrow Out(B_{\infty})$, which uses a version of $B_{\infty}$
in which the points to braid are the vertices of an (extended) binary tree. This
idea is further developed for $N$, where the group $A_N$ is defined in a more direct way (\cite{ka3}). In \S 3, we shall give a transparent geometric
description for the middle term of the extension $(NC)_T$ as well. This makes
use of the surface ${\cal S}_{\infty}$ which is the complement of a Cantor set in the
sphere $S^2$. We first prove (Theorem \ref{MC}) that $T$ embeds
into the mapping class group of this surface. The new geometric version of
$(NC)_T$ is then described in terms of mapping class groups (Definition-Proposition \ref{A}).\\

 We now give a statement of our main result (Theorem \ref{main}), which also concerns the Burau representation of the braid groups (see \cite{tu} for an
 up-to-date overview including recent results on the faithfulness question):\\

{\it Main Theorem -- The Burau representation of the Artin braid groups
  (depending on a parameter ${\bf t}\in\C^*$) extends
  to the mapping class group $A_T$. More precisely, there exists a Hilbert
  space ${\cal H}$ and a representation $\rho^{\bf t}: A_T\rightarrow GL({\cal
  H})$ (in the group of
  bounded operators) such that $B_{\infty}$ is represented in the subgroup
  ${\mathfrak T}$ of determinant operators. This gives a morphism of
  non-commutative extensions\\

 \setlength{\unitlength}{0.9cm}

\begin{picture}(10,2) 
\multiput(4,2)(1.5,0){2}{\vector(1,0){0.5}}   
\put(3.5,1.9){1} \put(4.7,1.9){$B_{\infty}$} \put(6.3,1.9){${\cal A}_T$}
\multiput(7,2)(2,0){2}{\vector(1,0){1}}   

\put(1.5,1.9){$(NC)_T$}

\put(8.4,1.9){$T$} 
\put(10.2,1.9){1}

\put(6.5,1.4){\vector(0,-1){.6}}  
\multiput(4,0.4)(1.2,0){2}{\vector(1,0){0.5}}   
\put(3.5,0.3){1} \put(4.8,0.3){${\mathfrak T}$} 
\put(6,0.3){$GL({\cal H})$}
\multiput(7.4,0.4)(2.2,0){2}{\vector(1,0){0.4}}   
\put(8.1,0.3){$\frac{GL({\cal H})}{\mathfrak T}$}
\put(10.2,0.3){1}

\put(1.5,0.3){$(NC)_{ps}$}
\put(4.9,1.4){\vector(0,-1){.6}}

\put(8.5,1.4){\vector(0,-1){.6}}
\end{picture}

inducing a morphism of central extensions\\

\setlength{\unitlength}{0.9cm}

\begin{picture}(10,2) 
\multiput(3.8,2)(1.5,0){2}{\vector(1,0){0.5}}   
\put(3.5,1.9){1} \put(4.7,1.9){$\Z$} \put(6.2,1.9){$\frac{{\cal A}_T}{[B_{\infty},B_{\infty}]}$}
\multiput(7.8,2)(1.5,0){2}{\vector(1,0){0.6}}   

\put(8.7,1.9){$T$}
\put(10.3,1.9){1}

\put(6.6,1.4){\vector(0,-1){.6}}  
\multiput(4,0.4)(1.2,0){2}{\vector(1,0){0.5}}   
\put(3.5,0.3){1} \put(4.6,0.3){$\C^*$} 
\put(6,0.3){$\frac{GL({\cal H})}{{\mathfrak T}_1}$}
\multiput(7.4,0.4)(2.2,0){2}{\vector(1,0){0.4}}   
\put(8.2,0.3){$\frac{GL({\cal H})}{\mathfrak T} $}
\put(10.3,0.3){1}

\put(4.8,1.4){\vector(0,-1){.6}}

\put(8.8,1.4){\vector(0,-1){.6}}
\end{picture}

where ${\mathfrak T}_1\subset {\mathfrak T}$ is the kernel of the determinant
morphism.}\\ 

The Pressley-Segal extension $(NC)_{PS}$ is a pull-back of the non-commutative
extension $(NC)_{ps}$ (see Definition-Proposition \ref{red}), which will be
called the {\it reduced Pressley-Segal extension} in our work. We also obtain a quite similar theorem linking the extension $(NC)_N$ to the reduced Pressley-Segal
extension $(NC)_{ps}$ (Theorem \ref{Neretin}).\\

We would like to mention that several interesting conections between the
Thompson (Neretin) groups and infinite dimensional groups and algebras have
been recently developed.\\
A. Reznikov \cite{re} (see also A. Navas \cite{na}) showed that the group $T$
has not the Kazhdan property. His argument uses in an essential way a
1-cocycle on the group $\mbox{Diff}^{1+\alpha}(S^1)$ (into which $T$ embeds,
see \cite{gh-se}) with values in the
Hilbert space of Hilbert-Schmidt operators, coming from the representation
of $\mbox{Diff}^{1+\alpha}(S^1)$ in $GL_{res}(L^2(S^1))$ when $\alpha >
\frac{1}{2}$.\\
D. Farley \cite{fa} later proved that the Thompson groups are a-T-menable (in
other words, have the Haagerup property). This led to a
proof of the Baum-Connes conjecture for these groups as well as to a new
argument for Reznikov's result.\\
Interesting representations of $N$ appear in \cite{ne} and \cite{ne'}, and embeddings of the Thompson groups in the Cuntz-Pimsner algebra have been
obtained by V. Nekrashevich \cite{nek} and J.-C. Birget \cite{bir}. Moreover, a $C^*$-algebra version of the class $\overline{gv}$ has been
studied by C. Oikonomides \cite{oi}.\\       

The article is organised as follows: Section 2 present in a convenient way the
background concerning the
Thompson and Neretin groups. Section 3 is mainly devoted to the construction
of the group $A_T$ as a mapping class group of a surface with infinitely many ends, but begins with a more elementary construction describing
$T$ itself as a mapping class group (Theorem \ref{MC}). In section 4, we recall the definition of
the Pressley-Segal extension, introduce its reduced version, and prove the main theorems (Theorems
\ref{Neretin}, \ref{main}).\\

{\bf Acknowledgements.} We warmly thank Louis Funar, Vaughan Jones and
Hitoshi Moriyoshi for very stimulating discussions at various stages of this work.

\section{Germ groups of inverse monoids: the Thompson and Nere\-tin groups}

{\bf 2.1. The inverse monoid of Fredholm tree automorphisms.} Let $T$ be a locally
finite tree, $\partial T$ be its boundary at infinity with
its usual topology, cf. e.g. \cite{ser}. A {\it cofinite domain} $D$ is the
complement of a finite subgraph of $T$. It has finitely many
components, which are subtrees of $T$.\\ 
A {\it Fredholm automorphism} of $T$, or {\it partial tree automorphism}
of $T$, is a bijection between two cofinite domains $D$ and $D'$ of $T$, which induces a tree
isomorphism on each connected component of $D$. The set
$Fred(T)$ of Fredholm automorphisms forms an inverse monoid (cf. \cite{pa}) with the obvious
structure: if $g$ and $h$ in $Fred(T)$ are defined on $D_g$ and $D_h$
respectively, $g. h$ is the Fredholm automorphism on the cofinite domain
$D_{g.h}:=h^{-1}(h(D_h)\cap D_g)$ defined by restriction of $g\circ h$. The
inverse $g^{-1}$ is defined on the cofinite domain $g(D_g)$.\\
For $g\in Fred(T)$, we call $D_g$ the source, and $g(D_g)$ the target.\\
The index of a Fredholm automorphism $g$ defined on the cofinite domain $D_g$ is the integer
$$ind\; g=card(Vert(T\setminus D_g))-card(Vert(T\setminus g(D_g))),$$
where $Vert(G)$ denotes the set of vertices of a graph $G$.\\
The subset $Fred^0(T)$ of Fredholm automorphisms with null index is an inverse
submonoid of $Fred(T)$.\\

\noindent{\bf 2.2. The group germification.} Since the source and target of a Fredholm
automorphism $g$ are cofinite, their boundaries at infinity coincide with that of the tree
$T$. Since $g$ acts by tree isomorphims outside a finite subgraph, it induces a homeomorphism $\partial g$ of $\partial T$. The
morphism of inverse monoids between $Fred(T)$ and the homeomorphism group of
the boundary $Homeo(\partial T)$ of $T$ (considered as an inverse monoid)
 $$g\in Fred(T)\mapsto \partial g \in Homeo(\partial T)$$
is the {\it group germification map}. Its image, denoted ${\cal G}(T)$,
is the {\it germ group} of the inverse monoid $Fred(T)$. Similarly, we shall
denote ${\cal G}^0(T)$ the image of $Fred^0(T)$.\\

\noindent{\bf 2.3. Fundamental examples: the Neretin and Thompson groups.}\\

\noindent{\bf 2.3.1.} Let ${\cal T}_2$ be the regular tree whose vertices are
all 3-valent (${\cal T}_2$ is called the dyadic tree, as its boundary may be
identified with the projective line on the field of dyadic numbers
$\Q_2$). Clearly, $Fred({\cal T}_2)=Fred^0({\cal T}_2)$, and the associated
germ group is {\it Neretin's spheromorphism group} $N$ defined in \cite{ne}
(cf. also \cite{ka1}, \cite{ka2}, \cite{ka3}).\\

\noindent{\bf 2.3.2.} Let ${\cal T}$ be the planar rooted tree, whose root is
2-valent, while the other vertices are 3-valent. Each vertex inherits a local
orientation from the orientation of the plane: the root being at the top of
the tree, each vertex $v$ has two descendants (or sons), the left one
$\alpha_l(v)$ and the right one $\alpha_r(v)$.\\ 
We shall say that $T$ is a {\it finite dyadic rooted subtree} of ${\cal T}$
(f.d.r.s.t. for short) if it is a finite subtree of ${\cal T}$, rooted in the
root of ${\cal T}$, such that its vertices are $3$-valent or $1$-valent
(except the root, which is $2$-valent). The $1$-valent vertices are called the
{\it leaves} of $T$. We denote the set of leaves of $T$ by ${\cal L}(T)$ and
say that the {\it type of $T$} is $k$ if it has $k$ leaves.\\    
The {\it canonical labelling} of a f.d.r.s.t whose type
is $k\in\N^*$ is the list of its leaves $v_0,\ldots,v_{k-1}$, enumerated
leftmost first, and reading from left to right (cf. Figure 1a).

\begin{de}[Thompson's group $T$, see also \cite{ca}]\label{tho}
Let $Fred^+({\cal T})$ be the inverse submonoid of $Fred({\cal
  T})=Fred^0({\cal  T})$ consisting of planar partial automorphisms, that is, $g$ belongs to $Fred^+({\cal T})$ if there exist two
  f.d.r.s.t. $T_0$ and $T_1$ of ${\cal T}$ such that
\begin{enumerate}
\item $D_g:=({\cal T}\setminus T_0)\cup {\cal L}(T_0)$ is the source of $g$;
\item  $g(({\cal T}\setminus T_0)\cup {\cal L}(T_0))=({\cal T}\setminus T_1)\cup
  {\cal L}(T_1)$;
\item $g$ induces a cyclic bijection ${\cal L}(T_0)\rightarrow {\cal L}(T_1)$ in the
  following sense: if $(v_0,\ldots,v_{k-1})$ and
  $(w_0,\ldots,w_{k-1})$ are the canonical labellings of the leaves of $T_0$
  and $T_1$ respectively, there exists a cyclic permutation $\sigma$ of the set
$\{0,\ldots,k-1\}$ identified with $\Z/k\Z$ such that $g(v_i)=w_{\sigma(i)}$
for $i=0,\ldots, k-1$. The cyclicity of $\sigma$ means that there exists some
$i_0\in \Z/k\Z$ such that for all $i$, $\sigma(i)=i+i_0$ mod $k$;
\item for each vertex $v$ of $D_g$,
  $g(\alpha_i(v))=\alpha_i(g(v))$, $i=l,r$.
\end{enumerate}
Thompson's group $T$ is then defined as the germ group of $Fred^+({\cal T})$. 
\end{de}

\noindent {\bf Symbols}. Each $g\in Fred^+({\cal T})$ may be uniquely represented by a symbol
$(T_1,T_0,\sigma)$, with the notations of the above definition. We denote by
$[T_1,T_0,\sigma]$ the element of $T$ which is the germ of $g$. We have the following composition rule in $T$:
$$[T_2, T_1,\sigma][T_1, T_0,\tau]=[T_2,T_0,\sigma\circ\tau].$$

\begin{example} Figure 1a represents a symbol coupled with a cyclic permutation
$\sigma$.
\end{example}

\begin{figure}
\begin{center}
\begin{picture}(0,0)%
\includegraphics{A.pstex}%
\end{picture}%
\setlength{\unitlength}{3315sp}%
\begingroup\makeatletter\ifx\SetFigFont\undefined%
\gdef\SetFigFont#1#2#3#4#5{%
  \reset@font\fontsize{#1}{#2pt}%
  \fontfamily{#3}\fontseries{#4}\fontshape{#5}%
  \selectfont}%
\fi\endgroup%
\begin{picture}(6669,2080)(274,-1361)
\put(3781,-151){\makebox(0,0)[lb]{\smash{\SetFigFont{8}{9.6}{\rmdefault}{\mddefault}{\updefault}$\sigma(0)=1$}}}
\put(3601,-106){\makebox(0,0)[lb]{\smash{\SetFigFont{10}{12.0}{\rmdefault}{\mddefault}{\updefault},}}}
\put(2026,-1321){\makebox(0,0)[lb]{\smash{\SetFigFont{8}{9.6}{\rmdefault}{\mddefault}{\updefault}Figure 1a}}}
\put(6166,-961){\makebox(0,0)[lb]{\smash{\SetFigFont{7}{8.4}{\rmdefault}{\mddefault}{\updefault}$\frac{1}{2}$   }}}
\put(6886,-961){\makebox(0,0)[lb]{\smash{\SetFigFont{7}{8.4}{\rmdefault}{\mddefault}{\updefault}$1$}}}
\put(5491,-961){\makebox(0,0)[lb]{\smash{\SetFigFont{7}{8.4}{\rmdefault}{\mddefault}{\updefault}$0$}}}
\put(5311,614){\makebox(0,0)[lb]{\smash{\SetFigFont{7}{8.4}{\rmdefault}{\mddefault}{\updefault}$1$}}}
\put(5311,-421){\makebox(0,0)[lb]{\smash{\SetFigFont{7}{8.4}{\rmdefault}{\mddefault}{\updefault}$\frac{1}{4}$}}}
\put(5311,-61){\makebox(0,0)[lb]{\smash{\SetFigFont{7}{8.4}{\rmdefault}{\mddefault}{\updefault}$\frac{1}{2}$   }}}
\put(6526,-961){\makebox(0,0)[lb]{\smash{\SetFigFont{7}{8.4}{\rmdefault}{\mddefault}{\updefault}$\frac{3}{4}$}}}
\put(5851,-1321){\makebox(0,0)[lb]{\smash{\SetFigFont{8}{9.6}{\rmdefault}{\mddefault}{\updefault}Figure 1b}}}
\put(2521,-106){\makebox(0,0)[lb]{\smash{\SetFigFont{10}{12.0}{\rmdefault}{\mddefault}{\updefault}$v_0$}}}
\put(2881,-466){\makebox(0,0)[lb]{\smash{\SetFigFont{10}{12.0}{\rmdefault}{\mddefault}{\updefault}$v_1$}}}
\put(3421,-466){\makebox(0,0)[lb]{\smash{\SetFigFont{10}{12.0}{\rmdefault}{\mddefault}{\updefault}$v_2$}}}
\put(1576,-106){\makebox(0,0)[lb]{\smash{\SetFigFont{10}{12.0}{\rmdefault}{\mddefault}{\updefault}$w_2$}}}
\put(2161,-106){\makebox(0,0)[lb]{\smash{\SetFigFont{10}{12.0}{\rmdefault}{\mddefault}{\updefault},}}}
\put(676,-466){\makebox(0,0)[lb]{\smash{\SetFigFont{10}{12.0}{\rmdefault}{\mddefault}{\updefault}$w_0$}}}
\put(1216,-466){\makebox(0,0)[lb]{\smash{\SetFigFont{10}{12.0}{\rmdefault}{\mddefault}{\updefault}$w_1$}}}
\end{picture}
\caption{An element of Thompson's group $T$} 
\end{center}
\end{figure}

\begin{de}[Thompson's group $F$, see also \cite{ca}]\label{de1}
Thompson's group $F$ is the germ group of the inverse submonoid of
$Fred^+({\cal T})$ whose elements are represented by symbols
$(T_1,T_0,\sigma)$, where $\sigma$ is the identity. It is a subgroup of $T$.
\end{de}

\noindent{\bf 2.4. Piecewise tree automorphisms.} If $U$ is a null index
  Fredholm operator of a Hilbert space ${\cal H}$, we may find a compact
  operator $K$ such that $U+K$ is an invertible operator. Analogous
  considerations in the present context lead to the following:

\begin{de-pr}
Let $Bij(T^0)$ be the group of bijections on the set of vertices $T^0$ of the
tree $T$. A piecewise tree automorphism of $T$ is a bijection $g\in Bij(T^0)$
which induces a Fredholm automorphism outside a finite subset of $T^0$. The set of piecewise
tree automorphisms of $T$ forms a subgroup of $Bij(T^0)$, denoted $PAut(T)$.\\
For each $g\in Fred^0(T)$, there exists $h\in PAut(T)$ inducing the same germ
as $g$ in the boundary $\partial T$.
\end{de-pr}

Denote by ${\mathfrak S}(T^0)$ the group of finitely supported permutations on
the set of vertices $T^0$. The fundamental objects resulting from the
preceding considerations are a ``non-commutative" extension of the germ group
${\cal  G}^0(T)$ together with an associated  central extension:

\begin{pr}\label{pr1}
There is a short exact sequence of groups
$$1\rightarrow {\mathfrak S}(T^0)\longrightarrow PAut(T)\longrightarrow {\cal
  G}^0(T)\rightarrow 1 $$
called the non-commutative extension of the germ group ${\cal
  G}^0(T)$. Denote by ${\mathfrak A}(T^0)$ the alternating subgroup of
${\mathfrak S}(T^0)$. After dividing by ${\mathfrak A}(T^0)$, the
non-commutative extension provides a central extension of ${\cal G}^0(T)$:

$$0\rightarrow \Z/2\Z\longrightarrow \widehat{{\cal  G}^0(T)}\longrightarrow {\cal
  G}^0(T)\rightarrow 1, $$
where $\widehat{{\cal  G}^0(T)}$ is the quotient group $PAut(T)/{\mathfrak
  A}(T^0)$.
\end{pr}

\begin{proof} The verification that ${\mathfrak  A}(T^0)$ is normal in $PAut(T)$ is
easy. Note that any extension by $\Z/2\Z$ is central. 
\end{proof}

\noindent{\bf 2.5. Discrete and combinatorial analogues of the Bott-Virasoro-Godbillon-Vey class.} The successful application of the
preceding proposition to the Neretin group $N={\cal G}^0({\cal T}_2)$ (cf. \S2.3.1) is the following

\begin{theo}(Combinatorial (or dyadic) analogue of the Bott-Virasoro-Godbillon-Vey class for the Neretin group, cf. \cite{ka3})\label{theo1}
The central extension of the Neretin spheromorphism group is non-trivial and
defines a class $\overline{GV}$ in $H^2(N,\Z/2\Z)$.
\end{theo}

However, the analogous construction for Thompson's group $T$ is trivial. Indeed,
at the price of slightly modifying the tree ${\cal T}$ to a tree $\tilde{\cal
  T}$ (cf. Remark \ref{rem2} below), we prove that 
the projection $PAut(\tilde{\cal T})\rightarrow {\cal G}(\tilde{\cal T})$
splits over $T$. The construction of the splitting $T\rightarrow PAut(\tilde{\cal T})$ relies on the interpretation of $T$ as a piecewise affine homeomorphism group of the circle that we now recall:\\

\begin{figure}
\begin{center}
\begin{picture}(0,0)%
\includegraphics{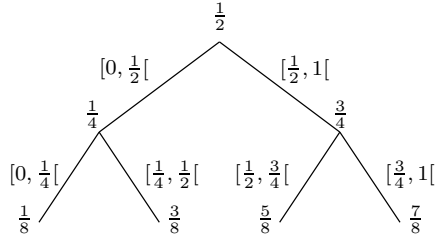}%
\end{picture}%
\setlength{\unitlength}{4144sp}%
\begingroup\makeatletter\ifx\SetFigFont\undefined%
\gdef\SetFigFont#1#2#3#4#5{%
  \reset@font\fontsize{#1}{#2pt}%
  \fontfamily{#3}\fontseries{#4}\fontshape{#5}%
  \selectfont}%
\fi\endgroup%
\begin{picture}(2385,1332)(1081,-703)
\put(3331,-421){\makebox(0,0)[lb]{\smash{\SetFigFont{8}{9.6}{\rmdefault}{\mddefault}{\updefault}$[\frac{3}{4},1[$}}}
\put(1531,-61){\makebox(0,0)[lb]{\smash{\SetFigFont{8}{9.6}{\rmdefault}{\mddefault}{\updefault}$\frac{1}{4}$}}}
\put(3016,-61){\makebox(0,0)[lb]{\smash{\SetFigFont{8}{9.6}{\rmdefault}{\mddefault}{\updefault}$\frac{3}{4}$}}}
\put(3466,-691){\makebox(0,0)[lb]{\smash{\SetFigFont{8}{9.6}{\rmdefault}{\mddefault}{\updefault}$\frac{7}{8}$}}}
\put(2296,524){\makebox(0,0)[lb]{\smash{\SetFigFont{8}{9.6}{\rmdefault}{\mddefault}{\updefault}$\frac{1}{2}$}}}
\put(1621,209){\makebox(0,0)[lb]{\smash{\SetFigFont{8}{9.6}{\rmdefault}{\mddefault}{\updefault}$[0,\frac{1}{2}[$}}}
\put(2431,-421){\makebox(0,0)[lb]{\smash{\SetFigFont{8}{9.6}{\rmdefault}{\mddefault}{\updefault}$[\frac{1}{2},\frac{3}{4}[$}}}
\put(1891,-421){\makebox(0,0)[lb]{\smash{\SetFigFont{8}{9.6}{\rmdefault}{\mddefault}{\updefault}$[\frac{1}{4},\frac{1}{2}[$}}}
\put(1081,-421){\makebox(0,0)[lb]{\smash{\SetFigFont{8}{9.6}{\rmdefault}{\mddefault}{\updefault}$[0,\frac{1}{4}[$}}}
\put(2701,209){\makebox(0,0)[lb]{\smash{\SetFigFont{8}{9.6}{\rmdefault}{\mddefault}{\updefault}$[\frac{1}{2},1[$}}}
\put(2566,-691){\makebox(0,0)[lb]{\smash{\SetFigFont{8}{9.6}{\rmdefault}{\mddefault}{\updefault}$\frac{5}{8}$}}}
\put(2026,-691){\makebox(0,0)[lb]{\smash{\SetFigFont{8}{9.6}{\rmdefault}{\mddefault}{\updefault}$\frac{3}{8}$}}}
\put(1126,-691){\makebox(0,0)[lb]{\smash{\SetFigFont{8}{9.6}{\rmdefault}{\mddefault}{\updefault}$\frac{1}{8}$}}}
\end{picture}

\caption{Labelling of the rooted dyadic tree by the dyadic rationals of
  $]0,1[$  }\label{fi2}
\end{center}
\end{figure}

\noindent{\bf Thompson's group $T$ as a group of piecewise dyadic affine homeomorphisms of the circle, cf. e.g. \cite{ca}.} Label inductively the edges of the tree ${\cal T}$ by all the dyadic intervals of $[0,1[$
in the following way: label the left descending edge of the root by
$[0,\frac{1}{2}[$, the right one by $[\frac{1}{2},1[$. If an edge
labelled by $[\frac{k}{2^n}, \frac{k+1}{2^n}[$
gives birth to two descending edges, the left one will be labelled by the first
half $[\frac{k}{2^n}, \frac{2k+1}{2^{n+1}}[$, and the right one by the second
half $[\frac{2k+1}{2^{n+1}}, \frac{k+1}{2^n}[$.\\
Label also the vertices by the dyadic rationals of $]0,1[$ inductively: the
root corresponds to $1/2$, and if an edge is labelled by
$[\frac{k}{2^n}, \frac{k+1}{2^n}[$, label its bottom vertex by the middle of
the interval, namely $\frac{2k+1}{2^{n+1}}$, cf. Figure \ref{fi2}.\\

\begin{de-pr}\label{de-pr0}
Let $g\in T$ be given by a symbol $(T_1,T_0,\sigma)$, with $\sigma$ a cyclic
bijection of $\Z/k\Z$ when the trees $T_0$ and $T_1$ have $k$ leaves. Denote
by $I_0,\ldots,I_{k-1}$ (resp. $J_0,\ldots,J_{k-1}$) the dyadic intervals corresponding to the terminal
edges of $T_0$ (resp. $T_1$). The unique piecewise affine map of $[0,1[$
applying affinely and increasingly $I_i$
onto $J_{\sigma(i)}$ for all $i=0,\ldots,k-1$, only depends on $g$, not on the
symbol. It induces an orientation-preserving homeomorphism $\phi_g$ of the
circle, viewed as $[0,1]/0\sim 1$. The correspondence $g_in T\mapsto \phi_g\in
Homeo^+(S^1)$ is a morphism, which embeds $T$ into the homeomorphism group of
the circle.  
\end{de-pr}

\begin{remark}\label{rem1} On each dyadic interval $I_i$ of the subdivision, $\phi_g$ is the
restriction of an affine map of the form $x\mapsto
2^{n_i}x+\frac{p_i}{2^{q_i}}$, $n_i\in \Z$, $p_i\in \Z$, $q_i\in \N$.
\end{remark} 

\begin{remark}\label{rem2} Add an edge $e_0$ to the tree ${\cal T}$, linked to the root of ${\cal T}$,
and label its terminal vertex by $0\sim 1$. Denote by $\tilde{\cal T}$ the
resulting tree. Thus, the vertices of $\tilde{\cal T}$ are in bijection
with the dyadic rationals of $[0,1]/0\sim 1$. Each $g\in T$ induces a
bijection $\hat{g}$ on the set of vertices $\tilde{\cal T}^0$ via the action
of $\phi_g$ on the set of dyadic rationals. The bijection $\hat{g}$ is a
piecewise tree automorphism, whose germ coincides with $g$, and the
correspondence $g\in T\rightarrow \hat{g}\in PAut(\tilde{\cal T})$ is a morphism.
\end{remark}
  
\begin{example}
Figure 1b represents $\phi_g$ for the element $g$ of $T$ defined on Figure 1a.\\
\end{example}

In \cite{gh-se}, the second cohomology group $H^2(T,\Z)$ is proved to be free
abelian on two generators $\chi$ and $\alpha$, where $\chi$ is the Euler class
of $T$, and $\alpha$ corresponds to the Godbillon-Vey class. In order to recall the formula of the
cocycle associated with $\alpha$, we need to define three functions on $S^1=[0,1]/0\sim 1$ associated with an element
$g\in T$: 

\begin{de}
Let $g$ be an element of $T$. For $x\in S^1=[0,1]/0\sim 1$, $g '_l(x)$ (resp. $g '_r(x)$) is the left
(resp. right) derivative number of the
affine bijection $\phi_g$ (it is an integral power of $2$), and $\Delta \log
_2 g'_r (x)$ is the integer $\log _2 g'_r (x)-\log _2 g'_l (x)$.
\end{de}

\begin{theo}(Discrete analogue of the Bott-Thurston cocycle for Thompson's group $T$,
  cf. \cite{gh-se})\\
The function $\overline{gv}: T\times T\rightarrow \Z$ defined by
  
$$\overline{gv}(g,h)=\sum_{x\in S^1} 
\begin{array}{|cc|}
\log _2 h'_r &\log _2 (g\circ h)'_r\\
\Delta \log _2 h'_r &\Delta \log _2 (g\circ h)'_r\\
 \end{array}
\;(x)$$ 
is a cocycle whose cohomology class equals $2\alpha$. Here $\begin{array}{|cc|}
a & b\\
c & d
 \end{array}$ is the determinant $ad-bc$.
\end{theo}

Note that the sum in the above formula is finite.\\

In the next section, we shall build geometrically a non-commutative
extension of Thompson's group $T$ (by a braid group) which abelianises to a
central extension by $\Z$ detecting the class $\alpha$. This
will be a braided version of the non-commutative extension of Neretin's group
(Proposition \ref{pr1} and \S 2.5).
   
\section{Geometric group extensions related to the genus zero infinite surface}

\noindent{\bf 3.1. The infinite surface.} Let $\bar{\cal S}_{\infty}$ be the genus zero infinite surface, constructed as an inductive
limit of finite subsurfaces $\bar{\cal S}_n$: $\bar{\cal S}_1$ is a compact cylinder, with two boundary
components, and $\bar{\cal S}_{n+1}$
is obtained from $\bar{\cal S}_n$ by gluing a copy of a ``pair of pants" (that is, a compact surface of genus zero with three boundary
components) along each boundary component of
$\bar{\cal S}_{n}$ (homeomorphic to a circle). It
follows that for each $n\geq 1$, $\bar{\cal S}_n$ is a $2^n$-holed sphere, and
$\bar{\cal S}_{\infty}=\displaystyle{\lim_{\stackrel{\rightarrow}{n}} \bar{\cal S}_n}$. The surface
$\bar{\cal S}_{\infty}$ is oriented, and a homeomorphism of the surface will always
be supposed to be orientation-preserving, unless the opposite is explicitly stated.\\

\begin{figure}\label{fi3}
\begin{center}
\begin{picture}(0,0)%
\includegraphics{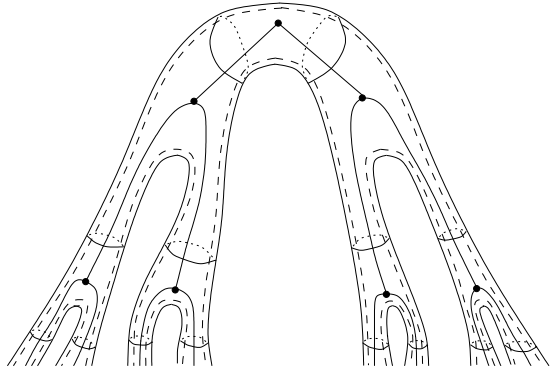}%
\end{picture}%
\setlength{\unitlength}{3315sp}%
\begingroup\makeatletter\ifx\SetFigFont\undefined%
\gdef\SetFigFont#1#2#3#4#5{%
  \reset@font\fontsize{#1}{#2pt}%
  \fontfamily{#3}\fontseries{#4}\fontshape{#5}%
  \selectfont}%
\fi\endgroup%
\begin{picture}(4074,2739)(214,-2098)
\end{picture}

\caption{Surface ${\cal S}_{\infty}$ and its tree}
\end{center}
\end{figure}

\noindent{\it Pants decomposition and rigid structure}: By this construction,
$\bar{\cal S}_{\infty}$ is naturally equipped with a pants
decomposition, which will be referred to in the sequel as the {\it canonical
decomposition}. We introduce a {\it rigid structure} on $\bar{\cal S}_{\infty}$,
consisting of three disjoint {\it seams} on each pair of pants, and two
disjoint seams on the cylinder $\bar{\cal S}_1$, as
indicated on Figure \ref{fi3}: a seam (represented by a dotted line on Figure
\ref{fi3}) is homeomorphic to a segment, and connects two boundary circles of the pair of pants or the
cylinder; each pair of boundary circles of a pair of pants is
connected by a unique seam; each seam extends continuously to the seams of the
adjacent pairs of pants or cylinder.\\
We fix also a marked point, or puncture, on each pair of pants as well as on the
cylinder. All marked pairs of pants are homeomorphic, by homeomorphisms which respect the seams and the
marked points. We call {\it rigid} such homeomorphisms.\\
We shall denote by ${\cal S}_{\infty}$ the surface $\bar{\cal S}_{\infty}$ decorated with the punctures.\\

\noindent{\bf 3.2. Associativity homeomorphisms.} We shall suppose that ${\cal S}_{\infty}$ is lying on an oriented
plane, in such a way that the punctures are drawn on the ``visible
 side" of the surface, and the seams separate the visible side from the hidden
 side. In other words, the visible and hidden sides are the connected
 components of the complement in ${\cal S}_{\infty}$ of the union of the
 seams.\\ 
The {\it tree of the surface $\bar{\cal S}_{\infty}$} (or ${\cal S}_{\infty}$)
 is the rooted dyadic
tree drawn on the visible side of the surface, whose vertices are the
 punctures of $\bar{\cal S}_{\infty}$, and edges are transverse to the circles
 of the canonical pants decomposition of $\bar{\cal S}_{\infty}$. It will be identified with the rooted dyadic planar
tree ${\cal T}$ of \S2.3.2.\\ 
By a {\it finite subsurface} of $\bar{\cal S}_{\infty}$ we shall always mean a connected
 finite union of pair of pants of the infinite surface, together with the
 cylinder. To each finite subsurface $S$ corresponds a f.d.r.s.t. $T_S$
 (cf. \S2.3.2), which is the subtree of ${\cal T}$ whose internal vertices are
 the punctures of $S$. We call it the {\it tree of the
   surface S}.\\
The canonical labelling of the leaves of the tree of a finite subsurface $S$
(cf. \S2.3.2) provides a canonical labelling of the set of
boundary components of $S$. We say that $S$ has {\it type $k\in \N^*$} if it has
$k$ boundary components. Equivalently, the type of $S$ is the type of its tree (cf. \S2.3.2).

\begin{de-pr}[associativity homeomorphisms]\label{de-pr1}
Let $S_0$ and $S_1$ be two finite subsurfaces of ${\cal S}_{\infty}$ having the
same type $k\in\N^*$, and $i$ an integer such
that $0\leq i\leq k-1$. Up to isotopy, there exists a unique homeomorphism
$\gamma_{1,0}^i:S_0\rightarrow S_1$ mapping the visible side of $S_0$ onto the
visible side of $S_1$, and the $0^{th}$ boundary component of $S_0$ onto the
$i^{th}$ boundary component of $S_1$.
\end{de-pr}

\begin{figure}
\begin{center}
\begin{picture}(0,0)%
\includegraphics{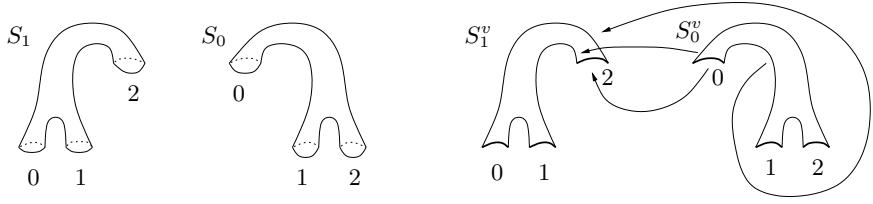}%
\end{picture}%
\setlength{\unitlength}{2901sp}%
\begingroup\makeatletter\ifx\SetFigFont\undefined%
\gdef\SetFigFont#1#2#3#4#5{%
  \reset@font\fontsize{#1}{#2pt}%
  \fontfamily{#3}\fontseries{#4}\fontshape{#5}%
  \selectfont}%
\fi\endgroup%
\begin{picture}(7381,1683)(901,-1033)
\put(3826,-916){\makebox(0,0)[lb]{\smash{\SetFigFont{9}{10.8}{\rmdefault}{\mddefault}{\updefault}$2$}}}
\put(1936,-196){\makebox(0,0)[lb]{\smash{\SetFigFont{9}{10.8}{\rmdefault}{\mddefault}{\updefault}$2$}}}
\put(1081,-916){\makebox(0,0)[lb]{\smash{\SetFigFont{9}{10.8}{\rmdefault}{\mddefault}{\updefault}$0$}}}
\put(1486,-916){\makebox(0,0)[lb]{\smash{\SetFigFont{9}{10.8}{\rmdefault}{\mddefault}{\updefault}$1$}}}
\put(3376,-916){\makebox(0,0)[lb]{\smash{\SetFigFont{9}{10.8}{\rmdefault}{\mddefault}{\updefault}$1$}}}
\put(2836,-196){\makebox(0,0)[lb]{\smash{\SetFigFont{9}{10.8}{\rmdefault}{\mddefault}{\updefault}$0$}}}
\put(5446,-871){\makebox(0,0)[lb]{\smash{\SetFigFont{9}{10.8}{\rmdefault}{\mddefault}{\updefault}$1$}}}
\put(5041,-871){\makebox(0,0)[lb]{\smash{\SetFigFont{9}{10.8}{\rmdefault}{\mddefault}{\updefault}$0$}}}
\put(7381,-826){\makebox(0,0)[lb]{\smash{\SetFigFont{9}{10.8}{\rmdefault}{\mddefault}{\updefault}$1$}}}
\put(7786,-826){\makebox(0,0)[lb]{\smash{\SetFigFont{9}{10.8}{\rmdefault}{\mddefault}{\updefault}$2$}}}
\put(6616,344){\makebox(0,0)[lb]{\smash{\SetFigFont{9}{10.8}{\rmdefault}{\mddefault}{\updefault}$S_0^v$}}}
\put(4816,299){\makebox(0,0)[lb]{\smash{\SetFigFont{9}{10.8}{\rmdefault}{\mddefault}{\updefault}$S_1^v$                 }}}
\put(901,299){\makebox(0,0)[lb]{\smash{\SetFigFont{9}{10.8}{\rmdefault}{\mddefault}{\updefault}$S_1$}}}
\put(2566,299){\makebox(0,0)[lb]{\smash{\SetFigFont{9}{10.8}{\rmdefault}{\mddefault}{\updefault}$S_0$}}}
\put(5986,-61){\makebox(0,0)[lb]{\smash{\SetFigFont{9}{10.8}{\rmdefault}{\mddefault}{\updefault}$2$}}}
\put(6931,-61){\makebox(0,0)[lb]{\smash{\SetFigFont{9}{10.8}{\rmdefault}{\mddefault}{\updefault}$0$}}}
\end{picture}

\caption{Construction of an associativity homeomorphism}\label{fi4}
\end{center}
\end{figure}

\begin{proof} Cutting $S_0$ along the seams, one gets two $2k$-gons, which are the
visible side $S_0^v$ and the hidden side $S_0^h$. Among the sides of the
$2k$-gon $S_0^v$, $k$ of them correspond to the boundary components of $S_0$,
from which they inherit the same labelling. Call them the distinguished sides
of $S_0^v$. Do the same with $S_1$. Since
$\gamma_{1,0}^i:S_0\rightarrow S_1$ must preserve the visible and hidden sides,
we need to define its restrictions $\gamma_{1,0}^{i,v}:S_0^v\rightarrow S_1^v$
and $\gamma_{1,0}^{i,h}:S_0^h\rightarrow S_1^h$. But clearly, there is a
unique $\gamma_{1,0}^{i,v}:S_0^v\rightarrow S_1^v$ which is
orientation-preserving, and maps the distinguished side $0$ of $S_0^v$ onto the
distinguished side $i$ of $S_1^v$. This forces the side $j$ to be mapped on
the side $j+i$ mod $k$ (cf. the example of Figure \ref{fi4}: $k=3$, $i=2$). In the same way, define
$\gamma_{1,0}^{i,h}:S_0^h\rightarrow S_1^h$  compatible with
$\gamma_{1,0}^{i,v}$ along the boundaries of the $2k$-gons. Thus,
$\gamma_{1,0}^{i,h}$ and $\gamma_{1,0}^{i,v}$ agree to induce the expected homeomorphism $\gamma_{1,0}^{i}:S_0\rightarrow S_1$. Its
unicity up to isotopy is clear. 
\end{proof}

\noindent{\bf 3.3. Thompson's group $T$ is a mapping class group.}

\begin{de}
A homeomorphism $\gamma$ of the surface $\bar{\cal S}_{\infty}$ is
asymptotically rigid if there exist two finite subsurfaces $S_0$ and $S_1$ of
$\bar{\cal S}_{\infty}$, having the same type, such that the restriction
$\gamma: \bar{\cal S}_{\infty}\setminus S_0\rightarrow \bar{\cal
  S}_{\infty}\setminus S_1$ is rigid, that is, maps rigidly a pair of pants
onto a pair of pants. 
\end{de}

\begin{theo}\label{MC}
Thompson's group $T$ embeds into the mapping class group of the surface
$\bar{\cal S}_{\infty}$ as the group of isotopy classes of those homeomorphisms representable
by asymptotically rigid homeomorphisms which preserve the visible side of $\bar{\cal
  S}_{\infty}$.
\end{theo}

\begin{proof}
Let $(T_1,T_0,\sigma)$ be a symbol defining an element
$g\in T$: $T_0$ and $T_1$ are $k$-ary trees, and
$\sigma\in\Z/k\Z$ prescribes a cyclic bijection from the set of leaves of $T_0$ to the set of leaves
of $T_1$. Let $S_0$ and $S_1$ be the finite subsurfaces of $\bar{\cal
  S}_{\infty}$ whose associated trees are respectively $T_0$ and $T_1$.  If $i=\sigma(0)$, let $\gamma^i_{1,0}:S_0\rightarrow S_1$ be the
associativity homeomorphism defined in \S3.2. Extend it to the unique
homeomorphism $\gamma_g:\bar{\cal S}_{\infty}\rightarrow
\bar{\cal S}_{\infty}$ which is rigid outside $S_0$. We claim that
$\gamma_g$  only depends on $g$, not on the choice of the symbol. Indeed,
since two symbols $(T_1,T_0,\sigma)$ and $(T'_1,T'_0,\sigma')$ defining the
same $g\in T$ always possess a common refining symbol $(T''_1,T''_0,\sigma'')$
(that is, such that the dyadic subdivision defined by $T''_0$ is a common
refinement of the subdivisions defined by $T_0$ and $T'_0$ respectively), it
is sufficient to suppose that $(T'_1,T'_0,\sigma')$ is a refinement of
$(T_1,T_0,\sigma)$, and by induction, that $(T'_1,T'_0,\sigma')$ is a simple
refinement of $(T_1,T_0,\sigma)$ (that is, $T'_0$ has just one more leaf than
$T_0$). Thus, the surface  $S'_0$ associated with $T'_0$ is the connected sum of $S_0$
with a pair of pants glued at some $j^{th}$ boundary component, while $S'_1$
associated with $T'_1$ is the connected sum of $S_1$
with a pair of pants glued at its $(j+i)^{th}$ (mod $k$) boundary
component. Denote by $\gamma_g':\bar{\cal S}_{\infty}\rightarrow\bar{\cal S}_{\infty}$ the
asymptotically rigid homeomorphism which rigidly extends the associativity
homeomorphism $\gamma^{i'}_{1,0}:S'_0\rightarrow S'_1$, where
$i'=\sigma'(0)$. At the price of replacing a homeomorphism by an isotopically
equivalent one, we may suppose that ${\gamma_g}_{|S'_0}=\gamma^{i'}_{1,0}$, hence
${\gamma_g'}_{|S'_0}={\gamma_g}_{|S'_0}$. Since by rigidity ${\gamma_g'}$ and
${\gamma_g}$ coincide outside $S'_0$, they coincide everywhere.\\
The unicity of ${\gamma_g}$ implies that the correspondence $g\in
T\mapsto {\gamma_g}\in MC(\bar{\cal S}_{\infty})$ (where
$MC(\bar{\cal S}_{\infty})$ denotes the mapping class group of the surface
$\bar{\cal S}_{\infty}$) is a morphism: if $g,g'$ belong to $T$, we may
represent $g$ by the symbol $(T',T,\sigma)$ and $g'$ by the symbol
$(T'',T',\sigma')$. It follows that $(T'',T,\sigma'\circ\sigma )$ is a symbol
for $g' g$. But clearly, $\gamma_{g'}\circ \gamma_g$ restricted to the finite
subsurface $S_{T}$ induces an associativity homeomorphism associated with the
symbol $(T'',T,\sigma'\circ\sigma )$, and it follows that $\gamma_{g'}\circ \gamma_g=\gamma_{g'.g}$. 
\end{proof}
 
\begin{remark}
The authors of \cite{fa-ga-ha} prove a similar theorem for Thompson's group
$F$.
\end{remark}

\noindent{\bf 3.4. Extension of Thompson's group $T$ by an infinite braid group.} \\

\begin{figure}
\begin{center}
\begin{picture}(0,0)%
\includegraphics{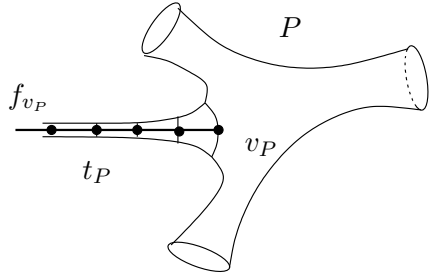}%
\end{picture}%
\setlength{\unitlength}{3729sp}%
\begingroup\makeatletter\ifx\SetFigFont\undefined%
\gdef\SetFigFont#1#2#3#4#5{%
  \reset@font\fontsize{#1}{#2pt}%
  \fontfamily{#3}\fontseries{#4}\fontshape{#5}%
  \selectfont}%
\fi\endgroup%
\begin{picture}(2795,1801)(676,-2045)
\put(1171,-1411){\makebox(0,0)[lb]{\smash{\SetFigFont{11}{13.2}{\rmdefault}{\mddefault}{\updefault}$t_P$}}}
\put(676,-916){\makebox(0,0)[lb]{\smash{\SetFigFont{11}{13.2}{\rmdefault}{\mddefault}{\updefault}$f_{v_P}$}}}
\put(2476,-466){\makebox(0,0)[lb]{\smash{\SetFigFont{11}{13.2}{\rmdefault}{\mddefault}{\updefault}$P$}}}
\put(2251,-1231){\makebox(0,0)[lb]{\smash{\SetFigFont{11}{13.2}{\rmdefault}{\mddefault}{\updefault}$v_P$}}}
\end{picture}

\caption{Pair of pants with tube and fibre}\label{fi5}
\end{center}
\end{figure}

\begin{figure}
\begin{center}
\begin{picture}(0,0)%
\includegraphics{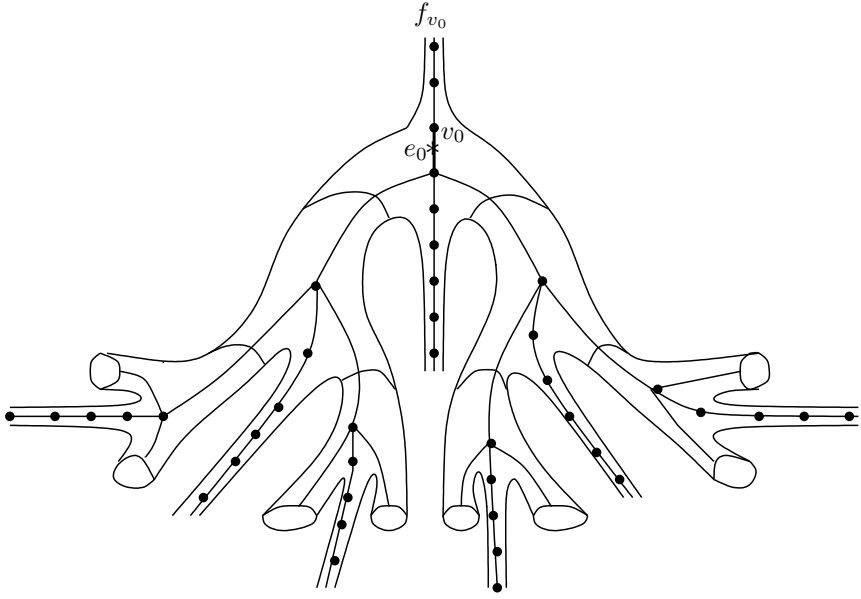}%
\end{picture}%
\setlength{\unitlength}{4972sp}%
\begingroup\makeatletter\ifx\SetFigFont\undefined%
\gdef\SetFigFont#1#2#3#4#5{%
  \reset@font\fontsize{#1}{#2pt}%
  \fontfamily{#3}\fontseries{#4}\fontshape{#5}%
  \selectfont}%
\fi\endgroup%
\begin{picture}(4270,2967)(2988,-2473)
\put(5176,-196){\makebox(0,0)[lb]{\smash{\SetFigFont{10}{12.0}{\rmdefault}{\mddefault}{\updefault}$v_0$}}}
\put(5041,389){\makebox(0,0)[lb]{\smash{\SetFigFont{10}{12.0}{\rmdefault}{\mddefault}{\updefault}$f_{v_0}$}}}
\put(4991,-286){\makebox(0,0)[lb]{\smash{\SetFigFont{10}{12.0}{\rmdefault}{\mddefault}{\updefault}$e_0$}}}
\put(5096,-286){\makebox(0,0)[lb]{\smash{\SetFigFont{12}{14.4}{\rmdefault}{\mddefault}{\updefault}$*$}}}
\end{picture}

\caption{Pants-with-tube decomposition of ${\cal S}_{\infty, t}$ with its tree
  ${\cal T}_t$}\label{fi6}
\end{center}
\end{figure}

\noindent{\bf 3.4.1. The infinite surface with tubes.} We want to build a new surface ${\cal S}_{\infty, t}$ by gluing some infinite
tubes on $\bar{\cal  S}_{\infty}$. Recall that $\bar{\cal  S}_{\infty}$ has a
visible and a hidden sides. In other words, there is an involutive
homeomorphism $j$ of $\bar{\cal  S}_{\infty}$, reversing the orientation,
stabilizing the circles of the canonical pants decomposition, whose set of
fixed points (the union of seams) bound two components of $\bar{\cal  S}_{\infty}$ (the visible and
hidden sides).\\
Let $P$ be a pair of pants of $\bar{\cal  S}_{\infty}$. Consider the seam
$s_P$ of $P$ connecting the two circles of its boundary which belong to the
boundary of the minimal finite subsurface containing both $P$ and the cylinder $\bar{\cal  S}_1$. Cut a small $j$-invariant disk on $P$
overlapping the seam $s_P$, and along the boundary of the resulting hole, glue an infinite
tubular surface. Cut also two small $j$-invariant disks on the cylinder $\bar{\cal
  S}_1$, each overlapping one of its two seams, and glue similarly two tubes along the
resulting holes.\\
Each tube is viewed as a connected countable union of compact cylinders. We add a
puncture on each boundary component of those cylinders, in such a way that the punctures are
alined along the tube (see Figure \ref{fi5}). If $P$ is a pair of pants of the surface
$\bar{\cal S}_{\infty}$ (resp. the cylinder $\bar{\cal S}_1$), denote by $t_P$
(resp. $t_P$ and $t_0$) the tubular surface (resp. surfaces) glued on $P$. The
puncture (resp. punctures) on the basis of $t_P$ (resp. $t_P$ and $t_0$) is
(resp. are) denoted by $v_P$ (resp. $v_P$ and $v_0$). We call the line passing by all the punctures the {\it fibre of the
  tube}, and denote it by $f_{v_P}$ ($f_{v_0}$). We denote by
${\cal S}_{\infty, t}$ the infinite surface with tubes and punctures (see
Figure \ref{fi6}).\\
The circles bounding the pants or the subcylinders of the tubes are called
the {\it circles of the pants-with-tubes decomposition of ${\cal S}_{\infty,
    t}$}.\\
Finally we may extend $j$ to an involutive homeomorphism $j_t$ of ${\cal S}_{\infty,
  t}$ (reversing the orientation and stabilizing the circles of the
pants-with-tubes decomposition of ${\cal S}_{\infty, t}$) and define the {\it
  visible side} of ${\cal S}_{\infty, t}$ as one of the two components bounded
by the set of fixed points of $j_t$. In particular, each tube has a visible
side, mapped by $j_t$ onto its hidden side, and we assume that all the fibres
belong to the visible side.\\   

\noindent{\bf Tree of the surface ${\cal S}_{\infty, t}$.} We draw the tree
$\tilde{\cal T}$ (see Remark \ref{rem2}) on the
visible side of ${\cal S}_{\infty, t}$, with vertices the punctures lying on
the bases of the tubes of the surface. The vertices of the added edge $e_0$ defined
in Remark \ref{rem2} are the punctures of the bases of the two tubes glued on
$\bar{\cal S}_1$. Thus, $v_0$ is the 1-valent vertex of $e_0$ in $\tilde{\cal T}$.\\ 
The union of $\tilde{\cal T}$ with all the fibres of the tubes is a tree ${\cal T}_{t}$
embedded in the visible side of the surface,  whose
vertices are the punctures of ${\cal S}_{\infty, t}$ (see Figure \ref{fi6}). We call ${\cal T}_{t}$ the
tree of the surface ${\cal S}_{\infty, t}$ and denote by $d$ its combinatorial metric.\\

\noindent{\bf Paths $\delta_v$}. Choose a base point $*$ on the surface ${\cal
  S}_{\infty, t}$, say, at the middle of the edge $e_0$. For
  each puncture $v$ of ${\cal S}_{\infty, t}$, that is, each vertex of
the tree ${\cal T}_{t}$, denote by $\delta_v$ the geodesic path from $*$ to
  $v$ contained in ${\cal T}_{t}$ (for its metric $d$).\\

\noindent{\bf 3.4.2. Asymptotically rigid homeomorphisms of the surface ${\cal
    S}_{\infty, t}$.}\\

By a {\it finite subsurface of ${\cal  S}_{\infty, t}$} we shall mean a
connected compact subsurface of ${\cal  S}_{\infty, t}$
containing the cylinder ${\cal S}_1$ and bounded by circles of the pants-with-tube
decomposition. The {\it bi-type} of a finite subsurface is the couple
of integers $(k,l)$, where $k$ is the number of circles of the boundary
components of $S$ which are boundaries of pants, and $l$ the number of
punctures on $S$. We check easily that $k$ is also the number of boundary
components which are circles of tubes, so that $S$ has exactly $2k$ boundary components.\\
If $S$ is a connected subsurface of ${\cal  S}_{\infty, t}$ (finite or not),
the {\it tree of $S$} is the trace of the tree ${\cal T}_t$ on $S$.\\
A homeomorphism $a$ of ${\cal  S}_{\infty, t}$ which permutes the punctures of
${\cal  S}_{\infty, t}$ is called {\it asymptotically
  rigid} if it verifies the following condition: there exist two finite subsurfaces $S_0$ and $S_1$ with the same
bi-type, such that $a$ restricts to each connected component $C$ of
${\cal  S}_{\infty, t}\setminus S_0$ to a {\it rigid} homeomorphism on its
image, that is, a homeomorphism mapping the visible side onto the visible
side, circles (of the canonical decomposition) onto circles, and the tree of
$C$ onto the tree of its image.\\
 
\noindent{\bf Germ of an asymptotically rigid homeomorphism}. It follows from the definition that an asymptotically rigid homeomorphism $a$
induces a piecewise tree automorphism of ${\cal T}_t$ (cf. \S2.4). This in
turn induces a homeomorphism of the boundary of
the tree ${\cal T}_t$, called the {\it
  germ} of $a$, and denoted by $\partial a$.\\

\noindent{\bf 3.4.3. The embedding $T\hookrightarrow Homeo(\partial {\cal
    T}_t)$}. Let $g$ be an element of $T$; as explained in section 2, Remark
\ref{rem2}, $g$ acts
on the set of vertices of $\tilde{\cal T}$, identified with the set of dyadic
rationals of $[0,1]/0\sim 1$. For any vertex $v$ of $\tilde{\cal T}$ labelled by the
dyadic rational $x$, denote the integer
$\Delta \log _2 g'_r(x)$ simply by $g''(v)$. Note that for almost all $x\in
[0,1]/0\sim 1$, $\Delta \log _2 g'_r(x)=0$, and we may check  that $\sum_x \Delta \log _2 g'_r(x)=0$, i.e. $\displaystyle{\sum_{v\in \tilde{\cal  T}^0} g''(v)=0}$. 

\begin{pr}
There is an embedding $T \hookrightarrow
Homeo(\partial {\cal T}_t )$, $g\mapsto \partial_{{\cal T}_t} g$, such that
each $\partial_{{\cal T}_t} g$ is induced by a Fredholm automorphism with null
index of ${\cal T}_t$.
\end{pr}

\begin{proof}  
If $g$ belongs to $T$, it is possible to extend the bijection $\hat{g}$ induced by $g$ on the set of vertices of
$\tilde{\cal T}$ (cf. Remark \ref{rem2}) to a Fredholm tree automorphism $\hat{g}_{{\cal T}_t}$ of ${\cal T}_t$  in the following
way:\\
 Let $N=\max_{v\in \tilde{\cal T}^0}
|g''(v)| +1$, and for each vertex $v$ of $\tilde{\cal T}$, denote by $f^{\geq N}_v$
the subfibre of $f_v$ whose vertices $s$ verify $d(v,s)\geq N$. 
\begin{enumerate}
\item If
$g''(v)=0$, define $\hat{g}_{{\cal T}_t}$ on $f_v$ as the unique isometric bijection from
$f_v$ to $f_{\hat{g}(v)}$;
\item if $g''(v)\not=0$, define $\hat{g}_{{\cal T}_t}$ on $f^{\geq N}_v$
as the unique isometric bijection from $f^{\geq N}_v$ to $f^{\geq  N+g''(v)}_{\hat{g}(v)}$.
\end{enumerate}
 Thus defined, $\hat{g}_{{\cal T}_t}$ is a Fredholm automorphism of
${\cal T}_t$, and its index is null, thanks to the identity $\sum_{v\in \tilde{\cal      T}^0} g''(v)=0$. The germ of $\hat{g}_{{\cal T}_t}$,
acting on the boundary of ${\cal T}_t$, only depends on $g$, and is
denoted $\partial_{{\cal T}_t} g$. The correspondence $$g\in T \mapsto
\partial_{{\cal T}_t} g \in Homeo(\partial {\cal T}_t )$$ is easily seen to be a group
morphism, thanks to the chain rule $(g\circ h)''(v)= g''(h(v))+ h''(v)$.
\end{proof}

\begin{remark}\label{cycl}
If $T_0$ and $T_1$ are finite subtrees of ${\cal T}_t$ such that
$\hat{g}_{{\cal T}_t}$ restricts to a bijection ${\cal T}_t\setminus
T_0\rightarrow {\cal T}_t\setminus T_1$, then $\hat{g}_{{\cal T}_t}$ induces a
{\it cyclic} bijection ${\cal L}(T_0)\rightarrow {\cal L}(T_1)$ between the
set of leaves of $T_0$ and $T_1$; the notion of cyclic bijection is quite
similar to that given in Definition \ref{tho}, {\it 3.}, since we may view
${\cal T}_t$ (embedded in the visible side of the infinite surface) as a planar tree.
\end{remark}

\noindent{\bf 3.4.4. The mapping class group ${\cal A}_T$.} Denote by  $Homeo^+({\cal
  S}_{\infty,t})$ the group of orientation-preserving homeomorphisms of the surface
  ${\cal  S}_{\infty,t}$ which permute the punctures.  

\begin{de-pr}[Mapping class group ${\cal A}_T$]\label{A}
Let ${\cal A}_T$ be the set of isotopy classes of homeomorphisms in
$Homeo^+({\cal S}_{\infty,t})$ which:
\begin{enumerate}
\item preserve the visible side of ${\cal S}_{\infty, t}$;
\item are asymptotically rigid in the sense of \S3.4.2; 
\item induce  partial tree automorphisms of the tree ${\cal T}_t$,
  whose germs belong to $T\hookrightarrow Homeo(\partial {\cal T}_t)$, cf. \S3.4.3.
\end{enumerate}
The set ${\cal A}_T$ is a group, and there is a short exact sequence 
$$1\rightarrow B_{\infty}[{\cal T}_t]\longrightarrow {\cal A}_{T} \longrightarrow T\rightarrow
1,$$ where $B_{\infty}[{\cal T}_t]$ denotes the mapping class group generated by the
isotopy classes of half-twists along the edges of ${\cal T}_t$, that is, the
braid group on the (planar) set of punctures
of ${\cal S}_{\infty, t}$.
\end{de-pr}

\begin{remark}\label{pres}
The tree ${\cal  T}_t$ being planar, a presentation of $B_{\infty}[{\cal  T}_t]$ may
  be deduced from \cite{se}, with generators the
half-twists between all pairs of consecutive vertices of the tree ${\cal T}_t$, along
the edge which connects them.
\end{remark}

\begin{remark}
The group ${\cal A}_{T}$ is defined in \cite{gr-se} in a combinatorial way,
by construction of a morphism $T\rightarrow Out(B_{\infty}[{\cal T}_t])$. The
main theorem of \cite{gr-se} asserts that the group ${\cal A}_{F'}$ (where
$F'$ is the derived subgroup of $F$), which results from the restriction of
the extension to $F'$, is acyclic.
\end{remark}

\begin{proof}
The set ${\cal A}_T$ is clearly a group and the natural
map ${\cal A}_T\rightarrow T$ a morphism. It remains to prove its surjectivity:
let $g$ be in $T$, $\hat{g}_{{\cal T}_t}$ the Fredholm automorphism of ${\cal
  T}_t$ which induces $g$ in the boundary of ${\cal T}_t$. Since
$\hat{g}_{{\cal T}_t}$ has null index, we may find two finite subsurfaces
$S_0$ and $S_1$ with the same bi-type $(k,l)$, and a rigid homeomorphism
$g_{0,1}:{\cal  S}_{\infty,t}  \setminus S_0\rightarrow {\cal  S}_{\infty,t}\setminus S_1$
(thus, preserving the visible sides), which induces $\hat{g}_{{\cal T}_t}$. We need to
extend $g_{0,1}$ to a $S_0$, such that it maps its visible side onto the
visible side of $S_1$.\\
For $i=0$ or 1, denote by $T_i$ the tree of the surface $S_i$, and by ${\cal
  L}(T_i)$ its set of leaves. The $2k$ boundary circles of $S_i$ are labelled
by the leaves of $T_i$. The point is that $g_{0,1}$ maps the circles of
$\partial S_0$ onto the circles of $\partial S_1$ by the prescription of the
bijection ${\cal L}(T_0)\rightarrow {\cal L}(T_1)$ induced by $\hat{g}_{{\cal
    T}_t}$ (see Remark \ref{cycl}). Since the latter is cyclic, we may proceed
as in the proof of Proposition \ref{de-pr1} to construct the extension of $g_{0,1}$
respecting the visible sides. The isotopy class of the resulting global
homeomorphism of ${\cal  S}_{\infty,t}$ belongs to ${\cal A}_T$ and its germ
is $g\in T$.\\
An element of the kernel of the projection ${\cal A}_T\rightarrow T$ may be
represented by a homeomorphism supported in a finite subsurface, preserving
its visible side and fixing pointwise its boundary circles. Thus, we may
suppose that its compact support is contained in the visible side. Since the homeomorphism can only
permute the punctures of the support, its isotopy class may be
identified with an Artin
braid on the punctures of the support. It follows that the kernel is the union
(or inductive limit) of the Artin braid groups defined for any such compact support. 
\end{proof}

\section{Pressley-Segal type representations of the geometric group extensions}

\subsection{The Pressley-Segal extension}

\noindent {\bf 4.1.1. Restricted linear group of a polarised Hilbert space.} Let
${\cal H}$ be a polarised separable Hilbert space, that is, the orthogonal sum
of two isometric separable Hilbert spaces  ${\cal H}={\cal H}_+\oplus{\cal
  H}_- $. Let $J$ be the operator $J=id_{{\cal H}_+}\oplus -id_{{\cal
    H}_-}$. The restricted algebra of bounded operators of ${\cal H}$ is
$${\cal L}_{res}({\cal H})=\{A\in {\cal L}({\cal H})|\;\; [J,A]\in {\cal
  L}_2({\cal H})\}.$$
 Here  ${\cal L}({\cal H})$ is the algebra of bounded
operators of ${\cal H}$, ${\cal  L}_2({\cal H})$ the ideal of
Hilbert-Schmidt operators, endowed with the usual Hilbert-Schmidt norm $||.||_2$. The restricted algebra is a Banach algebra for the
norm $||A||=|||A|||+||[J,A]||_2$.

\begin{de}[cf. \cite{pr-se}]
The restricted linear group $GL_{res}({\cal H})$ is the group of units of the
algebra ${\cal L}_{res}({\cal H})$. It is a Lie-Banach group.
\end{de}
    
Let $A=\left( \begin{array}{cc}
a & b\\
c & d
\end{array}
\right)$  be the block decomposition of $A\in {\cal L}({\cal H})$ relative to the
direct sum ${\cal H}={\cal H}_+\oplus{\cal H}_- $. The operator $[J,A]$ is
Hilbert-Schmidt if and only if $b$ and $c$ are Hilbert-Schmidt operators. If $A$ belongs
to $GL_{res}({\cal H})$, the invertibility of $A$
implies that $a$ is Fredholm in ${\cal H}_+$, and has an index $ind(a)\in \Z$. It is
easy to check that $ind: GL_{res}({\cal H})\rightarrow \Z$, $A\mapsto ind(a)$
is a morphism. It can be proved that it induces an isomorphism $\pi_0
(GL_{res}({\cal H}))\cong \Z$. Denote by $GL_{res}^0({\cal H})$ the connected
component of the identity. It is known to be a perfect group
(cf. \cite{co-ka}, \S5.4).\\

\noindent {\bf 4.1.2. Pressley-Segal extension of the restricted linear group.}
Let ${\cal L}_1({\cal H}_+)$ denote the ideal of trace-class operators of
${\cal H}_+$. It is a Banach algebra for the norm
$||b||_1=Tr(\sqrt{b^*b})$, where $Tr$
denotes the trace form. We
say that an invertible operator $q$ of ${\cal H}_+$ has a determinant if
$q-id_ {{\cal H}_+}=t$ is trace-class. Its determinant is the complex number
$det(q)=\sum_{i=0}^{+\infty} Tr(\wedge^i t)$, where $\wedge^i t$ is the
operator of the Hilbert space $\wedge^i {\cal H}_+$ induced by $t$ (cf. \cite{si}).\\

Denote by ${\mathfrak T}$ the subgroup of $GL({\cal H}_+)$ consisting of operators
which have a determinant, and by ${\mathfrak T}_1$ the kernel of the morphism
$det: {\mathfrak T}\rightarrow \C^*$.\\
Let ${\mathfrak A}$ be the subalgebra of ${\cal L}_{res}({\cal H})\times
   {\cal L}({\cal H}_+)$ consisting of couples $(A,q)$ such that $a-q$ is a
   trace-class operator of ${\cal H}_+$. It is a Banach algebra for the norm
   $||(A,q)||=||A||+||a-q||_1$. Let ${\mathfrak E}$ be the subgroup of $GL_{res}({\cal H})\times GL({\cal H}_+)$
consisting of the units of ${\mathfrak A}$. If  $(A,q)$ belongs to
${\mathfrak E}$,  then $ind(a)=ind(q+(a-q))=ind(q)=0$, so that $A$ belongs to $GL_{res}^0({\cal H})$, and there is a short exact sequence
$$(NC)_{PS}\hspace{1cm} 1\rightarrow {\mathfrak T}\stackrel{i}{\longrightarrow} {\mathfrak E}\stackrel{p}{\longrightarrow}
GL_{res}^0({\cal H})\rightarrow 1.$$
Here $p(A,q)=A$, and $i(q)=(id_{\cal H}, q)$.

\begin{remark}
The construction of $(NC)_{PS}$ can be extended to a short exact sequence $1\rightarrow
{\mathfrak T} \longrightarrow \tilde{\mathfrak E}\longrightarrow
GL_{res}({\cal H})\rightarrow 1$ (cf. \cite{pr-se}). Note, however, that if a
{\it simple} group is represented in $GL_{res}({\cal H})$, its image is contained in
$GL_{res}^0({\cal H})$. 
\end{remark}  

\begin{de}[cf. \cite{pr-se}]
The (non-central) extension  $1\rightarrow {\mathfrak T}\longrightarrow {\mathfrak E} \longrightarrow
GL_{res}^0({\cal H})\rightarrow 1$ will be called Pressley-Segal's
extension. It induces the central extension 
 
$$1\rightarrow \frac{{\mathfrak T}}{{\mathfrak T}_1 }\cong \C^* \longrightarrow
\frac{{\mathfrak E}}{{\mathfrak T}_1} \longrightarrow GL_{res}^0({\cal H})\rightarrow
1.$$
The corresponding cohomology class in
$H^2(GL_{res}^0({\cal H});\C^*)$ is denoted by $PS$, and called the
Pressley-Segal class of the restricted linear group.  
\end{de}

\noindent{\bf 4.1.3. Godbillon-Vey class and Pressley-Segal class.} Denote by
$\mbox{Diff}(S^1)$ the group of orientation-preserving diffeomorphisms of the
circle, and by $L^2(S^1)$ the Hilbert space of complex valued functions on the
circle which are square integrable. Then $L^2(S^1)$ admits a polarisation as
$L^2(S^1)={\cal H}_+\oplus {\cal H}_-$, where ${\cal H}_+$ is densely
generated by $\{e^{in\theta},\;n\in \Z, n\geq 0\}$ and  ${\cal H}_-$ by
$\{e^{in\theta},\;n\in \Z, n\leq -1\}$.\\
Consider the
series of representations, parametrised by $s\in\R$, $\pi_s:
\mbox{Diff}(S^1)\rightarrow GL^0_{res}(L^2(S^1))$,
$\pi_s(\phi)(f)=f\circ\psi. {\psi'}^s$, where $\psi=\phi^{-1}$. When $s=\frac{1}{2}$, the
representation is unitary. From computations of \cite{moriy}, the
pull-back of the Pressley-Segal class by $\pi_s$ is
$$\pi_s^*(PS)= e^{2i\pi(\frac{6s^2-6s+1}{6} gv -\frac{1}{12}\chi)},$$
where $gv\in H^2_{cont}(\mbox{Diff}(S^1);\R)$ is the Bott-Virasoro-Godbillon-Vey
class in the continuous cohomology of $\mbox{Diff}(S^1)$, and $\chi\in
H^2_{cont}(\mbox{Diff}(S^1);\R)$ is the continuous Euler class of
$\mbox{Diff}(S^1)$, that is, the smooth version of the discrete class of the extension $0\rightarrow \Z\longrightarrow \widetilde{\mbox{Diff}}(S^1)\longrightarrow \mbox{Diff}(S^1)\rightarrow1$, where $\widetilde{\mbox{Diff}}(S^1)$ is the universal cover of $\mbox{Diff}(S^1)$. Note that in the case $s=\frac{1}{2}$,
this result is consistent with \cite{pr-se}, Proposition 6.8.5.\\

\noindent{\bf 4.1.4. Reduced Pressley-Segal extension.}\\

The projection $p:A=\left( \begin{array}{cc}
a & b\\
c & d
\end{array}
\right)\in {\cal L}_{res}({\cal H})\mapsto  a\; mod\;{\cal L}_1({\cal H}_+)\in
\frac{{\cal L}({\cal H}_+) }{{\cal L}_1({\cal H}_+) }$ is a morphism of
algebras (essentially because the product of two Hilbert-Schmidt operators is
trace-class). It follows that the algebra $\mathfrak{A}$ is obtained by the
fibre-product

\begin{center}
\begin{picture}(0,0)%
\includegraphics{fibre.pstex}%
\end{picture}%
\setlength{\unitlength}{3315sp}%
\begingroup\makeatletter\ifx\SetFigFont\undefined%
\gdef\SetFigFont#1#2#3#4#5{%
  \reset@font\fontsize{#1}{#2pt}%
  \fontfamily{#3}\fontseries{#4}\fontshape{#5}%
  \selectfont}%
\fi\endgroup%
\begin{picture}(1812,1167)(361,-703)
\put(361,-691){\makebox(0,0)[lb]{\smash{\SetFigFont{8}{9.6}{\rmdefault}{\mddefault}{\updefault}${\cal L}({\cal H}_+)$}}}
\put(1981,-691){\makebox(0,0)[lb]{\smash{\SetFigFont{8}{9.6}{\rmdefault}{\mddefault}{\updefault}$\frac{{\cal L}({\cal H}_+)}{{\cal L}_1({\cal H}_+)}$}}}
\put(1891,299){\makebox(0,0)[lb]{\smash{\SetFigFont{8}{9.6}{\rmdefault}{\mddefault}{\updefault}${\cal L}_{res}({\cal H})$}}}
\put(1936,-196){\makebox(0,0)[lb]{\smash{\SetFigFont{8}{9.6}{\rmdefault}{\mddefault}{\updefault}$p$}}}
\put(586,299){\makebox(0,0)[lb]{\smash{\SetFigFont{8}{9.6}{\rmdefault}{\mddefault}{\updefault}$\mathfrak{A}$}}}
\end{picture}
\end{center}
  
Passing to the units, $p$ induces a group morphism $GL_{res}({\cal
  H})\rightarrow \left(\frac{{\cal L}({\cal H}_+)}{{\cal L}_1({\cal
      H}_+)}\right)^*$.
Denote by  $\left(\frac{{\cal L}({\cal H}_+)}{{\cal L}_1({\cal
  H}_+)}\right)_0^*$ the image of $GL^0_{res}({\cal  H})$ in
  $\left(\frac{{\cal L}({\cal H}_+)}{{\cal L}_1({\cal H}_+)}\right)^*$. 

\begin{de-pr}[Reduced Pressley-Segal extension]\label{red}
There is a morphism of non-commutative extensions

\begin{center}
\begin{picture}(0,0)%
\includegraphics{red.pstex}%
\end{picture}%
\setlength{\unitlength}{3315sp}%
\begingroup\makeatletter\ifx\SetFigFont\undefined%
\gdef\SetFigFont#1#2#3#4#5{%
  \reset@font\fontsize{#1}{#2pt}%
  \fontfamily{#3}\fontseries{#4}\fontshape{#5}%
  \selectfont}%
\fi\endgroup%
\begin{picture}(4005,844)(136,-290)
\put(2791,-241){\makebox(0,0)[lb]{\smash{\SetFigFont{8}{9.6}{\rmdefault}{\mddefault}{\updefault}$\left(\frac{{\cal L}({\cal H}_+)}{{\cal L}_1({\cal H}_+)}\right)^*_0$}}}
\put(2881,389){\makebox(0,0)[lb]{\smash{\SetFigFont{8}{9.6}{\rmdefault}{\mddefault}{\updefault}$GL^0_{res}({\cal H})$}}}
\put(1621,-241){\makebox(0,0)[lb]{\smash{\SetFigFont{8}{9.6}{\rmdefault}{\mddefault}{\updefault}$GL({\cal H}_+)$}}}
\put(811,389){\makebox(0,0)[lb]{\smash{\SetFigFont{8}{9.6}{\rmdefault}{\mddefault}{\updefault}${\mathfrak T}$ }}}
\put(811,-241){\makebox(0,0)[lb]{\smash{\SetFigFont{8}{9.6}{\rmdefault}{\mddefault}{\updefault}${\mathfrak T}$ }}}
\put(1711,389){\makebox(0,0)[lb]{\smash{\SetFigFont{8}{9.6}{\rmdefault}{\mddefault}{\updefault}${\mathfrak E}$}}}
\put(4141,-241){\makebox(0,0)[lb]{\smash{\SetFigFont{8}{9.6}{\rmdefault}{\mddefault}{\updefault}$1$}}}
\put(4141,389){\makebox(0,0)[lb]{\smash{\SetFigFont{8}{9.6}{\rmdefault}{\mddefault}{\updefault}$1$}}}
\put(136,389){\makebox(0,0)[lb]{\smash{\SetFigFont{8}{9.6}{\rmdefault}{\mddefault}{\updefault}$1$}}}
\put(136,-241){\makebox(0,0)[lb]{\smash{\SetFigFont{8}{9.6}{\rmdefault}{\mddefault}{\updefault}$1$}}}
\put(2971,119){\makebox(0,0)[lb]{\smash{\SetFigFont{8}{9.6}{\rmdefault}{\mddefault}{\updefault}$p$}}}
\end{picture}

\end{center}
 which induces a morphism of central extensions

\begin{center}
\begin{picture}(0,0)%
\includegraphics{cent.pstex}%
\end{picture}%
\setlength{\unitlength}{3315sp}%
\begingroup\makeatletter\ifx\SetFigFont\undefined%
\gdef\SetFigFont#1#2#3#4#5{%
  \reset@font\fontsize{#1}{#2pt}%
  \fontfamily{#3}\fontseries{#4}\fontshape{#5}%
  \selectfont}%
\fi\endgroup%
\begin{picture}(4095,844)(46,-290)
\put(2791,-241){\makebox(0,0)[lb]{\smash{\SetFigFont{8}{9.6}{\rmdefault}{\mddefault}{\updefault}$\left(\frac{{\cal L}({\cal H}_+)}{{\cal L}_1({\cal H}_+)}\right)^*_0$}}}
\put(2881,389){\makebox(0,0)[lb]{\smash{\SetFigFont{8}{9.6}{\rmdefault}{\mddefault}{\updefault}$GL^0_{res}({\cal H})$}}}
\put(4141,-241){\makebox(0,0)[lb]{\smash{\SetFigFont{8}{9.6}{\rmdefault}{\mddefault}{\updefault}$1$}}}
\put(4141,389){\makebox(0,0)[lb]{\smash{\SetFigFont{8}{9.6}{\rmdefault}{\mddefault}{\updefault}$1$}}}
\put(2971,119){\makebox(0,0)[lb]{\smash{\SetFigFont{8}{9.6}{\rmdefault}{\mddefault}{\updefault}$p$}}}
\put(1756,389){\makebox(0,0)[lb]{\smash{\SetFigFont{8}{9.6}{\rmdefault}{\mddefault}{\updefault}$\frac{{\mathfrak E}}{{\mathfrak T}_1}$}}}
\put( 46,389){\makebox(0,0)[lb]{\smash{\SetFigFont{8}{9.6}{\rmdefault}{\mddefault}{\updefault}$1$}}}
\put( 46,-241){\makebox(0,0)[lb]{\smash{\SetFigFont{8}{9.6}{\rmdefault}{\mddefault}{\updefault}$1$}}}
\put(586,389){\makebox(0,0)[lb]{\smash{\SetFigFont{8}{9.6}{\rmdefault}{\mddefault}{\updefault}$\frac{{\mathfrak T}}{{\mathfrak T}_1}\cong \C^*$ }}}
\put(1666,-241){\makebox(0,0)[lb]{\smash{\SetFigFont{8}{9.6}{\rmdefault}{\mddefault}{\updefault}$\frac{GL({\cal H}_+)}{{\mathfrak T}_1}$}}}
\put(586,-241){\makebox(0,0)[lb]{\smash{\SetFigFont{8}{9.6}{\rmdefault}{\mddefault}{\updefault}$\frac{{\mathfrak T}}{{\mathfrak T}_1}\cong \C^*$ }}}
\end{picture}
\end{center}
We call the non-commutative extension of $\left(\frac{{\cal L}({\cal H}_+)}{{\cal
      L}_1({\cal H}_+)}\right)^*_0$ the {\it Reduced Pressley-Segal
      extension}. The cohomology class in $H^2(\left(\frac{{\cal  L}({\cal
      H}_+)}{{\cal      L}_1({\cal H}_+)}\right)^*_0;\C^*)$ of the associated
      central extension is denoted by $ps$, and called the {\it Reduced
      Pressley-Segal class}. It follows that the class $PS$ is the pull-back of the reduced
      class $ps$ by the projection $p$.  
\end{de-pr}

The proof is straightforward. Since the ideal ${\cal L}_1({\cal H}_+)$ is not
closed in ${\cal L}({\cal H}_+)$, there is no separate topology on the quotient
$\frac{{\cal L}({\cal H}_+)}{{\cal L}_1({\cal H}_+)}$. This forces us to treat
the group $\left(\frac{{\cal L}({\cal H}_+)}{{\cal      L}_1({\cal
      H}_+)}\right)^*_0$, from a cohomological point of view, as a discrete
group.

\subsection{Pressley-Segal type representations for the non-com\-mutative and
  central extensions of Neretin's group}

In the next theorem, we refer to the material introduced in \S2.3.1, \S2.4 and
\S2.5, Theorem \ref{theo1}. 

\begin{theo}\label{Neretin}
Let ${\cal H}$ denote the Hilbert space $\ell^2({\cal T}_2^0)$ on the set of
vertices of the dyadic complete tree. There exists a unitary representation
$\rho: PAut({\cal T}_2)\rightarrow GL({\cal H})$ which induces a morphism of non-commutative extensions\\ 

\setlength{\unitlength}{0.9cm}
\begin{picture}(10,2) 
\multiput(4,2)(2,0){2}{\vector(1,0){0.5}}   
\put(3.5,1.9){1} \put(4.7,1.9){$\mathfrak{S}({\cal T}_2^0)$} \put(6.6,1.9){$PAut({\cal T}_2)$}
\multiput(8.4,2)(1.9,0){2}{\vector(1,0){0.9}}   

\put(9.6,1.9){$N$}
\put(11.4,1.9){1}

\put(4,0.4){\vector(1,0){0.5}}   
\put(5.4,0.4){\vector(1,0){0.8}}
\put(3.5,0.3){1} \put(4.9,0.3){${\mathfrak T}$} \put(6.6,0.3){$GL({\cal H})$}
\multiput(8.3,0.4)(2.5,0){2}{\vector(1,0){0.5}}   

\put(9,0.3){$\left(\frac{{\cal L}({\cal H)}}{{\cal      L}_1({\cal
      H})}\right)^*_0$}
\put(11.4,0.3){1}

\put(5,1.4){\vector(0,-1){.6}}
\put(7.1,1.4){\vector(0,-1){.6}}  
\put(9.7,1.4){\vector(0,-1){.6}}

\end{picture}

This, in turn, induces a morphism of central extensions\\

\setlength{\unitlength}{0.9cm}

\begin{picture}(10,2) 
\multiput(3.8,2)(1.8,0){2}{\vector(1,0){0.4}}   
\put(3.5,1.9){1} \put(4.4,1.9){$\Z/2\Z$} \put(6.3,1.9){$\frac{PAut({\cal
      T}_2) }{\mathfrak{A}({\cal T}_2^0)}$}
\multiput(8,2)(1.6,0){2}{\vector(1,0){0.5}}   

\put(8.9,1.9){$N$}
\put(10.4,1.9){1}

\multiput(3.9,0.4)(1.3,0){2}{\vector(1,0){0.5}}   
\put(3.5,0.3){1} \put(4.6,0.3){$\C^*$} 
\put(6.1,0.3){$\frac{GL({\cal H})}{{\mathfrak T}_1}$}
\multiput(7.4,0.4)(2.3,0){2}{\vector(1,0){0.4}}   
\put(8.1,0.3){$\left(\frac{{\cal L}({\cal H})}{{\cal      L}_1({\cal
      H})}\right)^*_0$}
\put(10.4,0.3){1.}

\put(4.8,1.4){\vector(0,-1){.6}}
\put(6.7,1.4){\vector(0,-1){.6}}  
\put(9,1.4){\vector(0,-1){.6}}
\end{picture}
\end{theo}

\begin{proof} Let $\{x_v\}_{v\in {\cal T}_2^0}$ be the hilbertian basis of ${\cal
    H}$ on the set of vertices of the tree ${\cal T}_2$. If $g$ belongs to
    $PAut({\cal T}_2)$, define the operator $\rho(g)$ on the basis by
    $\rho(g)(x_v)=x_{g(v)}$. This clearly defines a unitary operator of ${\cal    H}$. If $g$ is finitely supported, that is, $g$ belongs to ${\mathfrak S}({\cal
    T}_2^0)$, then $\rho(g)$ belongs to an inductive limit of finite
    dimensional linear groups $GL(n,\C)$, and in particular, to the group
    ${\mathfrak T}$. The commutativity of the second diagram essentially
    relies on
    the remark that the signature of $g\in {\mathfrak S}({\cal T}_2^0)$
    coincides with the determinant of the operator $\rho(g)$. 
\end{proof}

\subsection{A Pressley-Segal type representation of the mapping class group
    ${\cal A}_T$ extending the Burau representation}

\noindent{\bf 4.2.1. The Magnus and Burau representations.}\\

{\it Geometric definition of the Burau representation:} Let $D_n$ be the unit disk with $n$ marked points $q_1,\ldots,q_n$. Define the
index morphism by $ind: \pi_1(D_n)\rightarrow \Z$, $\gamma_i\mapsto 1$, where
$\gamma_1,\ldots,\gamma_n$ is a basis of the free group $\pi_1(D_n)$ (each
$\gamma_i$ is the homotopy class of a loop encircling $q_i$, based at a
chosen point $*$ on the boundary of $D_n$). The kernel of this epimorphism is the fundamental group of a covering space
$\tilde{D}_n$ over $D_n$. The transformation group of the covering
$\tilde{D}_n\rightarrow D_n$ is isomorphic to $\Z$. If $t$ denotes a generator
of this transformation group, viewed as a multiplicative group, then
$H_1(\tilde{D}_n)$ is a free module of rank $n-1$ over
$\Z[t,t^{-1}]$. Choose a point $\tilde{*}$ in $\tilde{D}_n$ over $*$. Since each $b\in
Homeo^+(D_n,\partial D_n)$ -- the group of orientation-preserving
homeomorphisms fixing pointwise the boundary of the disk -- acts on
$\pi_1(D_n)$ preserving the index, it has a unique lift $\tilde{b}$ to $\tilde{D}_n$ that
fixes $\tilde{*}$.  Hence there is a group morphism $Homeo^+(D_n,\partial D_n)\rightarrow
GL_{n-1}(\Z[t,t^{-1}])$, $b\mapsto \rho(b)$, where $\rho(b)$
is the action induced by $\tilde{b}$ on the homology group
$H_1(\tilde{D}_n)$. Since $\rho(b)$ only depends on the isotopy class of $b$,
one gets an induced morphism
$$ B_n\rightarrow
GL_{n-1}(\Z[t,t^{-1}]),$$
where $B_n=\pi_0(Homeo^+(D_n,\partial D_n))$ is the Artin braid group.\\

{\it Algebraic definition of the Burau representation: using the Magnus representation
  of $Aut(F_n)^{ind}$ (cf. \cite{bi}):}\\
Let $F_n$ be the free group on $\gamma_1,\ldots,\gamma_n$, and $Aut(F_n)^{ind}$ be the subgroup
  of $Aut(F_n)$ consisting of index-preserving automorphisms of $F_n$. There
  exists a representation ${\cal M}$ of $Aut(F_n)^{ind}$ in the linear group
  $GL_n(\Z[t,t^{-1}])$, called the {\it Magnus representation} (cf. \cite{bi}). On
  the other hand, consider the natural embedding of the braid
  group $B_n$ in $Aut(F_n)^{ind}$, deduced from the action of
  $Homeo^+(D_n,\partial D_n)$ on $\pi_1(D_n)=F_n$. It happens that the {\it Burau representation} is algebraically defined as the composition of the
  natural morphism $B_n\rightarrow Aut(F_n)^{ind}$ with the Magnus representation. The resulting
  representation splits into a one-dimensional trivial representation and an
  $(n-1)$-dimensional irreducible part: the latter is the one we have geometrically described.\\
Explicitly, the action of $B_n$ on $\pi_1(D_n)=F_n$ is defined through the
  following formula (as usual, $\sigma_1,\ldots,\sigma_{n-1}$ are the standard
  generators of $B_n$):
$$(\sigma_i)_*: \gamma_i\mapsto \gamma_i \gamma_{i+1}\gamma_i^{-1},$$
                $$\gamma_{i+1}\mapsto \gamma_i$$
                $$\gamma_j\mapsto \gamma_j,\;\;j\not=i,i+1.$$
The Magnus representation of $Aut(F_n)^{ind}$ may be elegantly introduced using
R. H. Fox's free differential calculus. To avoid unuseful length, we content
ourselves to give the result of the composition of the Magnus representation with
the embedding $B_n\hookrightarrow Aut(F_n)^{ind}$: introduce a basis of the
free module of rank
$n$ over $\Z[t,t^{-1}]$, ${\bf x_1},\ldots,{\bf x_n}$, (which may
be thought of as the homology classes of the lifts in
$\tilde{D}_n$, based at $\tilde{*}$, of $\gamma_1,\ldots,\gamma_n$
respectively). It follows that

$${\cal M}((\sigma_i)_*)({\bf x_i})=(1-t){\bf
  x_i}+t{\bf x_{i+1}},$$
$${\cal M}((\sigma_i)_*)({\bf x_{i+1}})={\bf x_i},$$
$${\cal M}((\sigma_i)_*)({\bf x_j})={\bf x_j},\;\;j\not=i,i+1.$$ 

Finally, $\rho(\sigma_i)\in GL_{n-1}(\Z[t,t^{-1}])$ is obtained by restricting ${\cal
  M}((\sigma_i)_*)$ to the submodule generated by ${\bf x_1}-{\bf x_2},\ldots,{\bf x_{n-1}}-{\bf x_n}$.\\  

\noindent{\bf 4.2.2. Extension of the Burau representation to the infinite braid
  group $B_{\infty}$}.\\

 For each puncture $q$ on the surface ${\cal  S}_{\infty,t}$, we define
  a loop $\gamma_q$ based at $*$ in the following way:\\

\begin{figure}
\begin{center}
\begin{picture}(0,0)%
\includegraphics{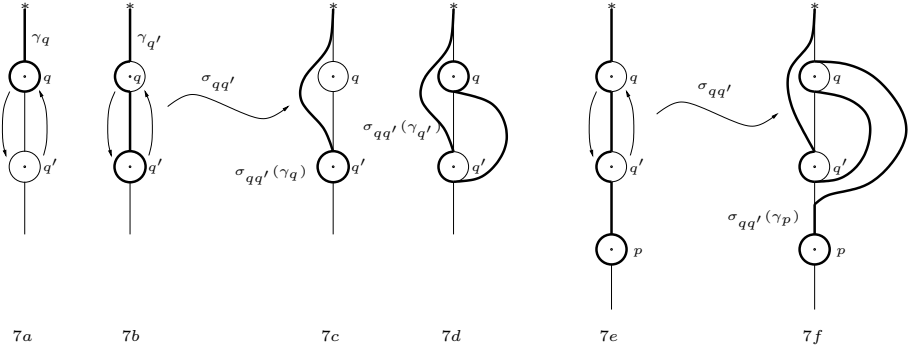}%
\end{picture}%
\setlength{\unitlength}{2072sp}%
\begingroup\makeatletter\ifx\SetFigFont\undefined%
\gdef\SetFigFont#1#2#3#4#5{%
  \reset@font\fontsize{#1}{#2pt}%
  \fontfamily{#3}\fontseries{#4}\fontshape{#5}%
  \selectfont}%
\fi\endgroup%
\begin{picture}(10859,4741)(172,-4007)
\put(541,209){\makebox(0,0)[lb]{\smash{\SetFigFont{5}{6.0}{\familydefault}{\mddefault}{\updefault}$\gamma_q$}}}
\put(1936,-1366){\makebox(0,0)[lb]{\smash{\SetFigFont{5}{6.0}{\familydefault}{\mddefault}{\updefault}$q'$}}}
\put(1756,-286){\makebox(0,0)[lb]{\smash{\SetFigFont{5}{6.0}{\familydefault}{\mddefault}{\updefault}$q$}}}
\put(5806,-286){\makebox(0,0)[lb]{\smash{\SetFigFont{5}{6.0}{\familydefault}{\mddefault}{\updefault}$q$}}}
\put(5806,-1366){\makebox(0,0)[lb]{\smash{\SetFigFont{5}{6.0}{\familydefault}{\mddefault}{\updefault}$q'$}}}
\put(8866,-1951){\makebox(0,0)[lb]{\smash{\SetFigFont{5}{6.0}{\familydefault}{\mddefault}{\updefault}$\sigma_{qq'}(\gamma_p)$}}}
\put(1621,-3391){\makebox(0,0)[lb]{\smash{\SetFigFont{6}{7.2}{\familydefault}{\mddefault}{\updefault}$7b$}}}
\put(316,-3391){\makebox(0,0)[lb]{\smash{\SetFigFont{6}{7.2}{\familydefault}{\mddefault}{\updefault}$7a$}}}
\put(4006,-3391){\makebox(0,0)[lb]{\smash{\SetFigFont{6}{7.2}{\familydefault}{\mddefault}{\updefault}$7c$}}}
\put(5446,-3391){\makebox(0,0)[lb]{\smash{\SetFigFont{6}{7.2}{\familydefault}{\mddefault}{\updefault}$7d$}}}
\put(7336,-3391){\makebox(0,0)[lb]{\smash{\SetFigFont{6}{7.2}{\familydefault}{\mddefault}{\updefault}$7e$}}}
\put(9766,-3391){\makebox(0,0)[lb]{\smash{\SetFigFont{6}{7.2}{\familydefault}{\mddefault}{\updefault}$7f$}}}

\put(2566,-286){\makebox(0,0)[lb]{\smash{\SetFigFont{5}{6.0}{\familydefault}{\mddefault}{\updefault}$\sigma_{qq'}$}}}
\put(4501,-871){\makebox(0,0)[lb]{\smash{\SetFigFont{5}{6.0}{\familydefault}{\mddefault}{\updefault}$\sigma_{qq'}(\gamma_{q'})$}}}
\put(8506,-376){\makebox(0,0)[lb]{\smash{\SetFigFont{5}{6.0}{\familydefault}{\mddefault}{\updefault}$\sigma_{qq'}$}}}
\put(676,-286){\makebox(0,0)[lb]{\smash{\SetFigFont{5}{6.0}{\familydefault}{\mddefault}{\updefault}$q$}}}
\put(676,-1366){\makebox(0,0)[lb]{\smash{\SetFigFont{5}{6.0}{\familydefault}{\mddefault}{\updefault}$q'$}}}
\put(4366,-286){\makebox(0,0)[lb]{\smash{\SetFigFont{5}{6.0}{\familydefault}{\mddefault}{\updefault}$q$}}}
\put(4366,-1366){\makebox(0,0)[lb]{\smash{\SetFigFont{5}{6.0}{\familydefault}{\mddefault}{\updefault}$q'$}}}
\put(7696,-286){\makebox(0,0)[lb]{\smash{\SetFigFont{5}{6.0}{\familydefault}{\mddefault}{\updefault}$q$}}}
\put(7696,-1366){\makebox(0,0)[lb]{\smash{\SetFigFont{5}{6.0}{\familydefault}{\mddefault}{\updefault}$q'$}}}
\put(7741,-2356){\makebox(0,0)[lb]{\smash{\SetFigFont{5}{6.0}{\familydefault}{\mddefault}{\updefault}$p$}}}
\put(406,569){\makebox(0,0)[lb]{\smash{\SetFigFont{6}{7.2}{\familydefault}{\mddefault}{\updefault}$*$}}}
\put(1666,569){\makebox(0,0)[lb]{\smash{\SetFigFont{6}{7.2}{\familydefault}{\mddefault}{\updefault}$*$}}}
\put(4096,569){\makebox(0,0)[lb]{\smash{\SetFigFont{6}{7.2}{\familydefault}{\mddefault}{\updefault}$*$}}}
\put(5536,569){\makebox(0,0)[lb]{\smash{\SetFigFont{6}{7.2}{\familydefault}{\mddefault}{\updefault}$*$}}}
\put(7426,569){\makebox(0,0)[lb]{\smash{\SetFigFont{6}{7.2}{\familydefault}{\mddefault}{\updefault}$*$}}}
\put(10126,-286){\makebox(0,0)[lb]{\smash{\SetFigFont{5}{6.0}{\familydefault}{\mddefault}{\updefault}$q$}}}
\put(10126,-1366){\makebox(0,0)[lb]{\smash{\SetFigFont{5}{6.0}{\familydefault}{\mddefault}{\updefault}$q'$}}}
\put(10171,-2356){\makebox(0,0)[lb]{\smash{\SetFigFont{5}{6.0}{\familydefault}{\mddefault}{\updefault}$p$}}}
\put(9856,569){\makebox(0,0)[lb]{\smash{\SetFigFont{6}{7.2}{\familydefault}{\mddefault}{\updefault}$*$}}}
\put(2971,-1411){\makebox(0,0)[lb]{\smash{\SetFigFont{5}{6.0}{\familydefault}{\mddefault}{\updefault}$\sigma_{qq'}(\gamma_q)$}}}
\put(1801,209){\makebox(0,0)[lb]{\smash{\SetFigFont{5}{6.0}{\familydefault}{\mddefault}{\updefault}$\gamma_{q'}$}}}
\end{picture}
\caption{Half-twists}\label{fi7}
\end{center}
\end{figure}

\noindent-- draw a small loop $c_q$ surrounding the puncture $q$;\\
-- slightly deform the edge-path ${\delta}_q$ to a path $\hat{\delta}_q$, such
  that the trace of $\hat{\delta}_q$ contains no punctures of ${\cal
  S}_{\infty,t}$, and $\hat{\delta}_q$ connects $*$ to the base point of the
  loop $c_q$. More precisely, we choose $\hat{\delta}_q$ such that it avoids the
  punctures by passing on their right (see Figure \ref{fi7}b, where
  $\hat{\delta}_{q'}$ passes on the right of $q$), and say that
  $\hat{\delta}_q$ has the {\it regularity property}.\\
--define $\gamma_q$ as the loop $\hat{\delta}_q.c_q.\hat{\delta}_q^{-1}$ based
  at $*$. Our convention is to compose the paths from left to right, so that
  $\hat{\delta}_q.c_q.\hat{\delta}_q^{-1}$ means
  $(\hat{\delta}_q.c_q).\hat{\delta}_q^{-1}$.\\
Let $F_{\infty}$ be the subgroup of $\pi_1({\cal S}_{\infty,t},*)$ generated
by the homotopy classes of loops $\gamma_q$, denoted ${\bf \gamma_q}$. It is a
free group on the set of punctures.\\
Recall (see Remark \ref{pres}) that $B_{\infty}[{\cal  T}_t]$ is generated the
half-twists between the pairs of consecutive vertices of the tree ${\cal T}_t$, along
the edge which connects them. So, denote by $\sigma_{q,q'}$ the half-twist
between two such consecutive vertices. The automorphism of $\pi_1({\cal S}_{\infty,t},*)$ induced
  by $\sigma_{q,q'}$ restricts to an automorphism $(\sigma_{q,q'})_*$ of the
  free subgroup $F_{\infty}$, that we now describe:
\begin{itemize}
\item $(\sigma_{q,q'})_*({\bf \gamma_q})={\bf \gamma_{q'}}$ (Figure 7c),
\item $(\sigma_{q,q'})_*({\bf \gamma_{q'}})={\bf \gamma_{q'}}{\bf
    \gamma_q}{\bf \gamma_{q'}}^{-1}$ (Figure 7d),
\item if $p\not=q,q'$ and $\gamma_p$ does not intersect the edge $qq'$,
  $(\sigma_{q,q'})_*({\bf \gamma_{p}})={\bf \gamma_{p}}$,
\item if $p\not=q,q'$ but $\gamma_p$ does intersect the edge $qq'$,
  $(\sigma_{q,q'})_*({\bf \gamma_{p}})={\bf \gamma_{q'}}{\bf
    \gamma_{q}}^{-1}{\bf \gamma_{p}}{\bf \gamma_{q}}{\bf \gamma_{q'}}^{-1}$
  (Figures 7e, 7f).
\end{itemize}

From the preceding formula we immediately deduce the

\begin{lem}
Define the index morphism $ind: F_{\infty}\rightarrow \Z$ by $ind ({\bf
  \gamma_{q}})=1$ for each puncture $q$ of ${\cal S}_{\infty,t}$. Each
  $\sigma\in B_{\infty }[{\cal T}_t]$ induces an automorphism $\sigma_*$ of
$F_{\infty}$ which is index-preserving.
\end{lem}

\begin{pr}
Let ${\cal H}$ be the Hilbert space $\ell^2({\cal T}_t^0)$ on the set
of punctures of the surface ${\cal S}_{\infty, t}$, that is, the vertices of the tree
${\cal T}_t$. Let ${\mathfrak T}$ be the group of determinant-operators of ${\cal H}$. For each complex
number ${\bf t}\in\C^*$, the Burau representation of the Artin braid groups extends
to a representation $\rho^{\bf t}_{\infty}:B_{\infty}[{\cal T}_t]\rightarrow {\mathfrak T}$. The composition with the determinant morphism induces an embedding
$H_1(B_{\infty}[{\cal T}_t])\hookrightarrow \C^*$ if and only if ${\bf t}$ is not a root
of unity.
\end{pr}

\begin{proof} Let $\{{\bf x_q}\}_q$ be a basis of the free $\Z[t,t^{-1}]$-module
on the set of punctures of ${\cal S}_{\infty,t}$. Each $\sigma\in
B_{\infty}[{\cal T}_t]$ induces $\sigma_*$ in $Aut(F_{\infty})^{ind}$ (the subgroup of
index-preserving automorphisms of $F_{\infty}$). Since there is no difficulty to define
the Magnus representation ${\cal M}_{\infty}$ in the infinite case, we define
$\rho_{\infty}(\sigma)$ in $Aut(\Z[t,t^{-1}]^{({\cal T}_t^0)})$ by
  $\rho_{\infty}(\sigma)={\cal M}_{\infty}(\sigma_*)$, and we get the formula: 
\begin{itemize}
\item if $p=q$, $\rho_{\infty}(\sigma)({\bf x_{q}})={\bf x_{q'}},$
\item if $p=q'$, $\rho_{\infty}(\sigma_)({\bf x_{q'}})=(1-t){\bf x_{q'}} +t{\bf x_{q}}$,
\item if $\gamma_p$, $p\not=q,q'$, intersects the edge between $q$ and $q'$,
  $\rho_{\infty}(\sigma)({\bf x_p})=(1-t)({\bf x_{q'}} -{\bf x_{q}}) +{\bf x_p} $,
\item in the other cases, $\rho_{\infty}(\sigma)({\bf x_p})={\bf x_p}$.
\end{itemize}
It follows that $\rho_{\infty}(\sigma)$ differs from the identity by a finite
rank operator. Note however that $\rho_{\infty}(\sigma)$ does not belong to
an inductive limit
$GL_{\infty}(\Z[t,t^{-1}])=\displaystyle{\lim_{\stackrel{\longrightarrow}{n}} GL_n (\Z[t,t^{-1}])}$,
  because an edge between $q$ and $q'$ intersects infinitely many loops
  $\gamma_p$.\\
By evaluation of $t$ on an invertible scalar ${\bf t}\in\C^*$, we get a representation in the
Hilbert space ${\cal H}$, $\rho^{\bf t}_{\infty}:B_{\infty}[{\cal T}_t]\rightarrow \mathfrak{T}$. Since $det\rho^{\bf t}_{\infty}(\sigma)=-{\bf t}$, $det\rho^{\bf
  t}_{\infty}$ induces the morphism $H_1(B_{\infty}[{\cal T}_t])\cong \Z\rightarrow \C^*$,
$n\mapsto (-{\bf t})^n$. 
\end{proof}      
       
\noindent{\bf 4.2.3. Pressley-Segal type representation of the mapping class
  group ${\cal A}_T$}. We are now ready to prove the main result of our paper:

\begin{theo}\label{main}
For each ${\bf t}\in\C^*$, the Burau representation $\rho^{\bf t}_{\infty}:
B_{\infty}:=B_{\infty}[{\cal T}_t]\rightarrow {\mathfrak T}$ extends to a
representation $\rho^{\bf t}_{\infty}$ of the mapping class group ${\cal A}_T$
in the Hilbert space ${\cal H}=\ell^2({\cal T}_t^0 )$ on the set
of punctures of the surface ${\cal S}_{\infty,t}$. There
is a morphism of non-commutative extensions\\
\setlength{\unitlength}{0.9cm}

\begin{picture}(10,2) 
\multiput(4,2)(1.5,0){2}{\vector(1,0){0.5}}   
\put(3.5,1.9){1} \put(4.7,1.9){$B_{\infty}$} \put(6.3,1.9){${\cal A}_T$}
\multiput(7,2)(2,0){2}{\vector(1,0){1}}   

\put(8.4,1.9){$T$} 
\put(10.2,1.9){1}

\put(6.5,1.4){\vector(0,-1){.6}}  
\multiput(4,0.4)(1.2,0){2}{\vector(1,0){0.5}}   
\put(3.5,0.3){1} \put(4.7,0.3){${\mathfrak T}$} 
\put(6,0.3){$GL({\cal H})$}
\multiput(7.4,0.4)(2.2,0){2}{\vector(1,0){0.4}}   
\put(7.9,0.3){$\left(\frac{{\cal L}({\cal H})}{{\cal L}_1({\cal      H})}\right)^*_0$}
\put(10.2,0.3){1}

\put(4.9,1.4){\vector(0,-1){.6}}

\put(8.5,1.4){\vector(0,-1){.6}}
\end{picture}

inducing a morphism of central extensions\\

\setlength{\unitlength}{0.9cm}

\begin{picture}(10,2) 
\multiput(3.8,2)(2,0){2}{\vector(1,0){0.4}}   
\put(3.5,1.9){1} \put(4.3,1.9){$H_1(B_{\infty})$} \put(6.3,1.9){$\frac{{\cal A}_T}{[B_{\infty},B_{\infty}]}$}
\multiput(7.8,2)(1.5,0){2}{\vector(1,0){0.9}}   

\put(8.9,1.9){$T$}
\put(10.3,1.9){1}

\put(6.8,1.4){\vector(0,-1){.6}}  
\multiput(3.9,0.4)(1.4,0){2}{\vector(1,0){0.5}}   
\put(3.5,0.3){1} \put(4.7,0.3){$\C^*$} 
\put(6.1,0.3){$\frac{GL({\cal H})}{{\mathfrak T}_1}$}
\multiput(7.4,0.4)(2.4,0){2}{\vector(1,0){0.4}}   
\put(8.1,0.3){$\left(\frac{{\cal L}({\cal H})}{{\cal L}_1({\cal      H})}\right)^*_0 $}
\put(10.3,0.3){1}

\put(4.9,1.4){\vector(0,-1){.6}}

\put(9,1.4){\vector(0,-1){.6}}
\end{picture}

where ${\mathfrak T}_1\subset {\mathfrak T}$ is the kernel of the determinant
morphism.\\
The vertical arrows are all injective whenever ${\bf t}\in\C^*$ is not a root
of unity.
\end{theo}

\begin{proof} It is subdivided into several lemmas.

\begin{lem}
There exists a representation ${\cal A}_T\rightarrow Aut(F_{\infty})$.
\end{lem}

\begin{proof} Let ${\cal A}_1$ be the group of isotopy classes of homeomorphisms of
${\cal S}_{\infty,t}$ which fix the base point, and satisfy the conditions
1--3 of Definition-Proposition \ref{A}. The isotopies are assumed to fix the base
point. Thus, there is a well-defined representation ${\cal A}_1\rightarrow
Aut(F_{\infty})$. One has an exact sequence $1\rightarrow
B_{\infty,1}\rightarrow {\cal A}_1\rightarrow T\rightarrow 1$, where
$B_{\infty,1}$ is generated by $B_{\infty}$ and a pure braid $\tau$, which can
be chosen as follows: let $v_*$ be one of the vertices of the distinguished
edge $e_0$, and $\sigma_{*v_*}$ be the half-twist along the half-edge which joins
$*$ to $v_*$. Define $\tau$ as the isotopy class of $\sigma_{*v_*}^2$ in
${\cal A}_1$. Next, consider the commutative diagram of short exact sequences:
\begin{center}
\begin{picture}(0,0)%
\includegraphics{K.pstex}%
\end{picture}%
\setlength{\unitlength}{2486sp}%
\begingroup\makeatletter\ifx\SetFigFont\undefined%
\gdef\SetFigFont#1#2#3#4#5{%
  \reset@font\fontsize{#1}{#2pt}%
  \fontfamily{#3}\fontseries{#4}\fontshape{#5}%
  \selectfont}%
\fi\endgroup%
\begin{picture}(3465,2685)(631,-2176)
\put(2431,-1681){\makebox(0,0)[lb]{\smash{\SetFigFont{8}{9.6}{\rmdefault}{\mddefault}{\updefault}${\cal A}_T$}}}
\put(2431,-916){\makebox(0,0)[lb]{\smash{\SetFigFont{8}{9.6}{\rmdefault}{\mddefault}{\updefault}${\cal A}_1$}}}
\put(2476,-151){\makebox(0,0)[lb]{\smash{\SetFigFont{8}{9.6}{\rmdefault}{\mddefault}{\updefault}$K$}}}
\put(3376,-916){\makebox(0,0)[lb]{\smash{\SetFigFont{8}{9.6}{\rmdefault}{\mddefault}{\updefault}$T$}}}
\put(3331,-1681){\makebox(0,0)[lb]{\smash{\SetFigFont{8}{9.6}{\rmdefault}{\mddefault}{\updefault}$T$}}}
\put(1171,-916){\makebox(0,0)[lb]{\smash{\SetFigFont{8}{9.6}{\rmdefault}{\mddefault}{\updefault}$B_{\infty, 1}$}}}
\put(1351,-151){\makebox(0,0)[lb]{\smash{\SetFigFont{8}{9.6}{\rmdefault}{\mddefault}{\updefault}$K$}}}
\put(631,-916){\makebox(0,0)[lb]{\smash{\SetFigFont{8}{9.6}{\rmdefault}{\mddefault}{\updefault}$1$}}}
\put(631,-1681){\makebox(0,0)[lb]{\smash{\SetFigFont{8}{9.6}{\rmdefault}{\mddefault}{\updefault}$1$}}}
\put(4096,-916){\makebox(0,0)[lb]{\smash{\SetFigFont{8}{9.6}{\rmdefault}{\mddefault}{\updefault}$1$}}}
\put(4096,-1681){\makebox(0,0)[lb]{\smash{\SetFigFont{8}{9.6}{\rmdefault}{\mddefault}{\updefault}$1$}}}
\put(2521,-2176){\makebox(0,0)[lb]{\smash{\SetFigFont{8}{9.6}{\rmdefault}{\mddefault}{\updefault}$1$}}}
\put(1396,-2176){\makebox(0,0)[lb]{\smash{\SetFigFont{8}{9.6}{\rmdefault}{\mddefault}{\updefault}$1$}}}
\put(1396,344){\makebox(0,0)[lb]{\smash{\SetFigFont{8}{9.6}{\rmdefault}{\mddefault}{\updefault}$1$}}}
\put(2521,344){\makebox(0,0)[lb]{\smash{\SetFigFont{8}{9.6}{\rmdefault}{\mddefault}{\updefault}$1$}}}
\put(1261,-1681){\makebox(0,0)[lb]{\smash{\SetFigFont{8}{9.6}{\rmdefault}{\mddefault}{\updefault}$B_{\infty}$}}}
\end{picture}
\end{center}

 The vertical arrows ${\cal A}_1\rightarrow {\cal A}$ and $B_{\infty,
 1}\rightarrow B_{\infty}$ are induced by
 forgetting the base point, and $K$, the kernel of those forgetting morphisms,
 is the normal subgroup of $B_{\infty,  1}$ generated by $\tau$. But it is
 easy to check that $\tau$ induces the identity in $Aut(F_{\infty})$. It
 follows that the representation ${\cal A}_1\rightarrow Aut(F_{\infty})$
 descends to ${\cal A}_1/K= {\cal A}_T\rightarrow
 Aut(F_{\infty})$. 
\end{proof}
      
Let $[a]$ be in ${\cal A}_T$, and choose any representative $a$ which fixes the
base point. There exist two finite subsurfaces $S_0$ and $S_1$
with the same bi-type, such that $a_{|{\cal S}_{\infty,t}\setminus S_0}:{\cal S}_{\infty,t}\setminus S_0\rightarrow {\cal S}_{\infty,t}\setminus S_1$ is rigid,
the action on ${\cal T}_t\cap ({\cal S}_{\infty,t}\setminus S_0)$ being prescribed by a partial tree
automorphism of ${\cal T}_t$. Label by $q_1,\ldots,q_n$ the punctures lying on
the boundary components of $S_0$, and $C_1,\ldots,C_n$ the corresponding
connected components of ${\cal S}_{\infty,t}\setminus S_0$ (note that $n=2k$
if $(k,l)$ is the bi-type of $S_0$).

\begin{lem}
For $i=1,\ldots,n$, there exists ${\bf \lambda_i}$ in $F_{\infty}$ such that
if $q$ is a puncture on the connected component $C_i$, then $a_*({\bf
  \gamma_q})= {\bf \lambda_i} {\bf  \gamma_{a(q)}}{\bf
  \lambda_i}^{-1}$. Furthemore, the homotopy class  ${\bf \lambda_i}$ belongs
to the free subgroup of $F_{\infty}$ generated by the ${\bf \gamma_q}$'s, for
$q\in {\cal T}_t^0\cap S_1$.  
\end{lem}

\begin{proof} Compare the homotopy classes, with fixed extremities, of
$a(\hat{\delta}_{q_i})$ and $\hat{\delta}_{a(q_i)}$. The loop $\lambda_i=a(\hat{\delta}_{q_i}).\hat{\delta}_{a(q_i)}^{-1}$ defines an
element ${\bf \lambda_i}$ of $F_{\infty}$.\\
For each $q$ on $C_i$, decompose the path $\hat{\delta}_q$ as
$\hat{\delta}_{q_i}.\hat{\delta}_{q_iq}$, where $\hat{\delta}_{q_iq}$ connects
$q_i$ to $q$. Thus, ${\bf \gamma_q}={\bf \hat{\delta}_{q_i}}{\bf
  \hat{\delta}_{q_iq}} {\bf c_q}{\bf  \hat{\delta}_{q_iq}}^{-1}{\bf
  \hat{\delta}_{q_i}}^{-1}$. Since the homeomorphism $a$ is rigid on $C_i$, it
preserves the regularity property of $\hat{\delta}_{q_iq}$ (see \S4.2.2), so
that $a(\hat{\delta}_{q_iq})= \hat{\delta}_{a(q_i)a(q)}$. It follows easily that
$a_*({\bf \gamma_q})= {\bf \lambda_i} {\bf  \gamma_{a(q)}}{\bf
  \lambda_i}^{-1}$.  
\end{proof}
           
\begin{lem} 
For each of the finitely many punctures $q$ of $S_0$, there exists
${\lambda_q}$ in $F_{\infty}$ such that $a_*({\bf \gamma_q})={\bf
  \lambda_q}{\bf \gamma_{a(q)}} {\bf \lambda_q} ^{-1}$.
\end{lem}

\begin{proof} The arguments are similar to those of the proof of the previous
  lemma. 
\end{proof}

It follows from the preceding three lemmas that each $a\in {\cal A}_T$ induces an index-preserving
automorphism of the free group $F_{\infty}$. Thus, we extend the Burau
representation to ${\cal A}_T$ by setting
$$\rho^{\bf t}_{\infty}(a)={\cal M}(a_*),$$
where ${\cal M}$ is, as before, the Magnus representation of $Aut(F_{\infty})^{ind}$. It remains to check that the
induced operator $\rho^{\bf t}_{\infty}(a)$ is indeed continuous:\\

For $i=1,\ldots,n$ and $\epsilon=0$ or $1$, introduce the subbasis ${\cal
  B}_i^{\epsilon}=\{{\bf x_q}, q\in {\cal T}_t^0\cap C_i^{\epsilon}\}$, where
  $C_i^0:=C_i$, and $C_i^1:=a(C_i)$, as well as the finite subbasis ${\cal
  B}_{n+1}^{\epsilon}=\{{\bf x_q}, q\in {\cal T}_t^0\cap S_{\epsilon}\}$.\\
If $\lambda_i=\gamma_{q_{j_1}}^{\epsilon_1}\ldots \gamma_{q_{j_r}}^{\epsilon_r}$ for
  some punctures $q_{j_1},\ldots q_{j_r}$ of $S_1$, and
  $\epsilon_1,\ldots,\epsilon_r\in\{-1,1\}$, then for any puncture $q$ lying on
  $C_i$:

$$\rho_{\infty}(a)({\bf x_q})=P_i(t,t^{-1})+\varepsilon_i t^{m_i} {\bf x_{a(q)}},$$
where $P_i(t,t^{-1})$ is a vector-valued polynomial in $t,t^{-1}$,
$\varepsilon_i\in\{-1,+1\}$, and $m_i\in\Z$, which only depend on $i$.\\
Similarly, for each of the finitely many $q$'s on $S_0$, there exists a vector-valued
polynomial $P_q(t,t^{-1})$,  an integer $m_q\in \Z$ and
$\varepsilon_q\in\{-1,+1\}$ such that
$$\rho_{\infty}(a)({\bf x_q})=P_q(t,t^{-1}) +\varepsilon_q t^m_q {\bf x_{a(q)}}.$$
Notice that the vectorial coefficients of the polynomials $P_i$'s and $P_q$'s belong
to the finite dimensional vector space spanned by the basis ${\cal
  B}_{n+1}^1$. It follows that the matrix of $\rho_{\infty}^{\bf t}(a)$ relative to the bases
${\cal B}^0=({\cal B}_1^0,\ldots,{\cal B}_{n+1}^0)$ and ${\cal B}^1=({\cal
  B}_1^1,\ldots,{\cal B}_{n+1}^1)$ may be written 

\begin{center}
\begin{picture}(0,0)%
\includegraphics{matrice.pstex}%
\end{picture}%
\setlength{\unitlength}{3522sp}%
\begingroup\makeatletter\ifx\SetFigFont\undefined%
\gdef\SetFigFont#1#2#3#4#5{%
  \reset@font\fontsize{#1}{#2pt}%
  \fontfamily{#3}\fontseries{#4}\fontshape{#5}%
  \selectfont}%
\fi\endgroup%
\begin{picture}(3522,1824)(1,-1423)
\put(2161,209){\makebox(0,0)[lb]{\smash{\SetFigFont{7}{8.4}{\rmdefault}{\mddefault}{\updefault}$0$}}}
\put(3241,209){\makebox(0,0)[lb]{\smash{\SetFigFont{7}{8.4}{\rmdefault}{\mddefault}{\updefault}$0$}}}
\put(2881,209){\makebox(0,0)[lb]{\smash{\SetFigFont{7}{8.4}{\rmdefault}{\mddefault}{\updefault}$0$}}}
\put(2161,-1276){\makebox(0,0)[lb]{\smash{\SetFigFont{9}{10.8}{\rmdefault}{\mddefault}{\updefault}$*$}}}
\put(2521,-1276){\makebox(0,0)[lb]{\smash{\SetFigFont{9}{10.8}{\rmdefault}{\mddefault}{\updefault}$*$}}}
\put(3241,-1276){\makebox(0,0)[lb]{\smash{\SetFigFont{9}{10.8}{\rmdefault}{\mddefault}{\updefault}$*$}}}
\put(2881,-1276){\makebox(0,0)[lb]{\smash{\SetFigFont{9}{10.8}{\rmdefault}{\mddefault}{\updefault}$*$}}}
\put(2161,-871){\makebox(0,0)[lb]{\smash{\SetFigFont{7}{8.4}{\rmdefault}{\mddefault}{\updefault}$0$}}}
\put(2521,-871){\makebox(0,0)[lb]{\smash{\SetFigFont{7}{8.4}{\rmdefault}{\mddefault}{\updefault}$0$}}}
\put(3241,-871){\makebox(0,0)[lb]{\smash{\SetFigFont{7}{8.4}{\rmdefault}{\mddefault}{\updefault}$0$}}}
\put(3241,-556){\makebox(0,0)[lb]{\smash{\SetFigFont{7}{8.4}{\rmdefault}{\mddefault}{\updefault}$0$}}}
\put(3241,-196){\makebox(0,0)[lb]{\smash{\SetFigFont{7}{8.4}{\rmdefault}{\mddefault}{\updefault}$0$}}}
\put(2881,-151){\makebox(0,0)[lb]{\smash{\SetFigFont{7}{8.4}{\rmdefault}{\mddefault}{\updefault}$0$}}}
\put(2521,-151){\makebox(0,0)[lb]{\smash{\SetFigFont{7}{8.4}{\rmdefault}{\mddefault}{\updefault}$0$}}}
\put(2881,-511){\makebox(0,0)[lb]{\smash{\SetFigFont{7}{8.4}{\rmdefault}{\mddefault}{\updefault}$0$}}}
\put(2521,-511){\makebox(0,0)[lb]{\smash{\SetFigFont{11}{13.2}{\rmdefault}{\mddefault}{\updefault}...}}}
\put(2161,-511){\makebox(0,0)[lb]{\smash{\SetFigFont{7}{8.4}{\rmdefault}{\mddefault}{\updefault}$0$}}}
\put(1711,209){\makebox(0,0)[lb]{\smash{\SetFigFont{5}{6.0}{\rmdefault}{\mddefault}{\updefault}$\pm t^{m_1}$}}}
\put(2071,-151){\makebox(0,0)[lb]{\smash{\SetFigFont{5}{6.0}{\rmdefault}{\mddefault}{\updefault}$\pm  t^{m_2}$}}}
\put(2791,-871){\makebox(0,0)[lb]{\smash{\SetFigFont{5}{6.0}{\rmdefault}{\mddefault}{\updefault}$\pm t^{m_n}$}}}
\put(2521,209){\makebox(0,0)[lb]{\smash{\SetFigFont{7}{8.4}{\rmdefault}{\mddefault}{\updefault}$0$}}}
\put(1846,-151){\makebox(0,0)[lb]{\smash{\SetFigFont{7}{8.4}{\rmdefault}{\mddefault}{\updefault}$0$}}}
\put(1846,-511){\makebox(0,0)[lb]{\smash{\SetFigFont{7}{8.4}{\rmdefault}{\mddefault}{\updefault}$0$}}}
\put(1846,-871){\makebox(0,0)[lb]{\smash{\SetFigFont{7}{8.4}{\rmdefault}{\mddefault}{\updefault}$0$}}}
\put(1846,-1276){\makebox(0,0)[lb]{\smash{\SetFigFont{9}{10.8}{\rmdefault}{\mddefault}{\updefault}$*$}}}
\put(  1,-286){\makebox(0,0)[lb]{\smash{\SetFigFont{7}{8.4}{\rmdefault}{\mddefault}{\updefault}$Mat_{({\cal B}_0,{\cal B}_1)}(\rho^{\bf t}_{\infty}(a))=$}}}
\end{picture}

\end{center}  
 
Each transformation  $\ell^2({\cal B}_i^0)\longrightarrow \ell^2( {\cal
  B}_i^1)$, ${\bf x_q}\mapsto \pm {\bf t}^{m_i} {\bf x_{a(q)}}$,
  $i=1,\ldots,n$, is obviously continuous. Since the $n+1$ blocks $[*]$ define
  finite rank, hence continuous  operators, it follows that $\rho^{\bf
  t}_{\infty}(a)$ is continuous.\\
The commutativity of the diagrams involved in the statement of Theorem \ref{main} is easy to
  prove. 
\end{proof}

\subsection{Concluding remarks}

It would be interesting to know if Theorem
\ref{main} has a version involving  the ``non-reduced" Pressley-Segal
extension. On the other hand,
we do not know if there is an extension of $\mbox{Diff}(S^1)$ which is the counterpart of our geometric extension of $T$ by $B_{\infty}$. What could be the kernel of it?\\
More generally, can one develop a representation theory
of $T$? Finally, we point out that the results of this paper do not
shed any light on Thompson's group $V$, acting on the Cantor set. However,
geometric group extensions of $V$ were recently studied by L. Funar and
the first author (\cite{fu-ka}), and by M. Brin and J. Meier.

\end{document}